\documentclass[12pt]{amsart}

\usepackage{times,epsf,amssymb,amsmath,hyperref, wrapfig, rlepsf, color}

\begin{document}

\newtheorem{thm}{Theorem}[section]
\newtheorem{lem}[thm]{Lemma}
\newtheorem{cor}[thm]{Corollary}
\newtheorem{conj}[thm]{Conjecture}

\theoremstyle{definition}
\newtheorem{defn}[thm]{\bf{Definition}}

\theoremstyle{remark}
\newtheorem{rmk}[thm]{Remark}

\def\square{\hfill${\vcenter{\vbox{\hrule height.4pt \hbox{\vrule width.4pt height7pt \kern7pt \vrule width.4pt} \hrule height.4pt}}}$}

\newcommand{\SI}{\partial_\infty ({\Bbb H}^2\times {\Bbb R})}
\newcommand{\Si}{S^1_{\infty}\times {\Bbb R}}
\newcommand{\si}{S^1_{\infty}}
\newcommand{\PI}{\partial_{\infty}}

\newcommand{\BH}{\Bbb H}
\newcommand{\BHH}{{\Bbb H}^2\times {\Bbb R}}
\newcommand{\BR}{\Bbb R}
\newcommand{\BC}{\Bbb C}
\newcommand{\BZ}{\Bbb Z}
\newcommand{\caps}{\overline{\BH^2}\times\{\pm\infty\}}

\newcommand{\e}{\epsilon}

\newcommand{\wh}{\widehat}
\newcommand{\wt}{\widetilde}

\newcommand{\A}{\mathcal{A}}
\newcommand{\C}{\mathcal{C}}
\newcommand{\p}{\mathcal{P}}
\newcommand{\R}{\mathcal{R}}
\newcommand{\B}{\mathcal{B}}
\newcommand{\h}{\mathcal{H}}
\newcommand{\T}{\mathfrak{T}}
\newcommand{\s}{\mathcal{S}}
\newcommand{\U}{\mathcal{U}}
\newcommand{\V}{\mathcal{V}}
\newcommand{\X}{\mathcal{X}}
\newcommand{\tc}{\mathcal{T}}
\newcommand{\Y}{\mathbf{Y}}

\newenvironment{pf}{{\it Proof:}\quad}{\square \vskip 12pt}

\title[Minimal Surfaces with Arbitrary Topology in $\BHH$]{Minimal Surfaces with Arbitrary Topology in $\BHH$}
\author{Baris Coskunuzer}
\address{UT Dallas, Dept. Math. Sciences, Richardson, TX 75080}
\email{coskunuz@utdallas.edu}
\thanks{The author is partially supported by Simons Collaboration Grant, and Royal Society Newton Mobility Grant.}

\maketitle

\begin{abstract} We show that any open orientable surface $S$ can be properly embedded in $\BHH$ as an area minimizing surface. 
\end{abstract}

\section{Introduction}

Minimal surfaces in $\BHH$ has been an attractive topic for the last two decades. After Nelli and Rosenberg's seminal results \cite{NR} on minimal surfaces in $\BHH$, the theory has been flourished very quickly with the substantial results on the existence, regularity, and other properties of minimal and CMC surfaces in $\BHH$, e.g. \cite{CR, CMT,HNST, KM, MMR,  MRR, NSST, RT, ST}.

In this paper, we are interested in the following question: {\em "What type of surfaces can be embedded into $\BHH$ as a complete minimal surface?"} 
Ros conjectured that any open orientable surface can be properly embedded in $\BHH$ as a minimal surface \cite{MR}. In this paper, we prove this conjecture.

\begin{thm} \label{main1} Any open orientable surface $S$ can be properly embedded in $\BHH$ as a complete area minimizing surface.
\end{thm}

In particular, this implies that any open orientable surface $S$ can be realized as an complete, embedded, minimal surface in $\BHH$. The key step is to show a vertical bridge principle for tall curves in $\Si$ (Section \ref{secbridge}). Then, by using the positive solutions of the asymptotic Plateau problem, we give a general construction to obtain complete, properly embedded minimal surfaces in $\BHH$ with arbitrary topology, i.e. any (finite or infinite) number of genus and ends.

The outline of the method is as follows: We start with a simple exhaustion of the open orientable surface $S$, i.e. $S_1\subset S_2 \ ...\ S_n \subset ..$ where $S=\bigcup_{n=1}^\infty S_n$. In particular, the surface $S$ is constructed by starting with a disk $D=S_1$, and by adding $1$-handles iteratively, i.e. $S_{n+1}-int(S_n)$ is either a pair of pants  or a cylinder with a handle (See Figure \ref{simpleexhaustion}). Hence after proving a bridge principle at infinity for $\BHH$ for vertical bridges in $\SI$ (Theorem \ref{bridge}), we started the construction with an area minimizing plane $\Sigma_1$ in $\BHH$. Then by following the iterative process dictated by the simple exhaustion, if $S_{n+1}$ is a pair of pants attached to $S_n$, then we attach one vertical bridge in $\SI$ to the corresponding component of $\PI \Sigma_n$. Similarly, if $S_{n+1}$ is a cylinder with a handle attached to $S_n$, then we attach two vertical bridges successively to $\PI \Sigma_n$ (See Figure \ref{hanger}) so that the number of boundary components of $\partial \Sigma_n$ and $\partial \Sigma_{n+1}$ are the same. By iterating this process, we inductively construct a properly embedded minimal surface $\Sigma$ in $\BHH$ with the same topological type of $S$.

The organization of the paper is as follows. In the next section, we give some definitions, and introduce the basic tools which we use in our construction. In Section \ref{secbridge}, we show the bridge principle at infinity in $\BHH$ for sufficiently long vertical bridges. In Section \ref{secsurfaceconstruction}, we prove the main result above. In Section \ref{secremarks}, we discuss generalization of our result to $H$-surfaces and finite total curvature case. We postpone some technical steps to the appendix at the end, where we also prove a generic uniqueness result for tall curves.

\subsection{Acknowledgements}
Part of this research was carried out at MIT and Max-Planck Institute during our visit. The author would like to thank them for their great hospitality.

\section{Preliminaries} \label{secprem}

In this section, we introduce our setup, and the basic tools which we use in our construction.

Throughout the paper, $\overline{\BHH}=\overline{\BH^2}\times\overline{\BR}=\BHH\cup\SI$ represents the natural product compactification of $\BHH$. In particular, $\SI=(\Si)\cup(\caps)$ represents the asymptotic boundary of $\BHH$. Also, we  call $\Si$ as the asymptotic cylinder, and $\caps$ as the caps at infinity.

\vspace{.2cm}

\noindent {\bf Convention:} Throughout the paper, by {\em curve}, we mean a finite collection of smooth Jordan curves unless otherwise stated.

\vspace{.2cm}

A curve $\Gamma$ in $\SI$ is \textit{finite} if $\Gamma \subset \Si$. If $\Gamma \cap \caps\neq \emptyset$, we call $\Gamma$ is {\em infinite}. Throughout the paper, all the curves in $\SI$ will be finite curves unless stated otherwise. For the asymptotic Plateau problem for infinite curves, see \cite{KM,Co3}.

\begin{defn} A compact surface with boundary $\Sigma$ is called {\em area minimizing surface} if $\Sigma$ has the smallest area among surfaces with the same boundary. A noncompact surface  is called {\em area minimizing surface} if any compact subsurface is an area minimizing surface.
\end{defn}

For our construction, one of our key ingredients is the solutions of the following Problem.

\vspace{.2cm}

\noindent {\bf Asymptotic Plateau Problem in $\BHH$:}

\vspace{.2cm}

\noindent {\em Let $\Gamma$ be a  collection of Jordan curves in $\Si$. Does there exist a complete, embedded minimal surface $\Sigma$ in $\BHH$ with $\PI\Sigma = \Gamma$?}

\vspace{.2cm}

Here, $\Sigma$ is an open, complete surface in $\BHH$, and $\PI\Sigma$ represents the asymptotic boundary of $\Sigma$ in $\SI$. Then,  $\overline{\Sigma}$ is the closure of $\Sigma$ in $\overline{\BHH}$, then $\PI \Sigma= \overline{\Sigma}\cap \SI$. Here, we stated the most general version of this problem. There are various results on this problem in the literature. For our construction, we need the positive solutions in a special case: Tall Curves (Lemma \ref{APP}).

\begin{defn}  \label{talldef} [Tall Curves] Consider the asymptotic cylinder $\Si$ with the coordinates $(\theta, t)$ where $\theta\in [0,2\pi)$ and $t\in \BR$. We call a rectangle $R=[\theta_1,\theta_2]\times[t_1,t_2]\subset \Si$  {\em tall rectangle} if $t_2-t_1>\pi$.
	
We call a finite collection of disjoint simple closed curves $\Gamma$ in $\Si$  {\em tall curve} if the region $\Gamma^c=\Si-\Gamma$ can be written as a union of open tall rectangles $R_i=(\theta^i_1,\theta^i_2)\times(t^i_1,t^i_2)$, i.e. $\Gamma^c=\bigcup_i R_i$.

We call a region $\Omega$ in $\Si$ a {\em tall region}, if $\Omega$ can be written as a union of tall rectangles, i.e. $\Omega=\bigcup_i R_i$ where  $R_i$ is a tall rectangle.
\end{defn}

Note that tall rectangles in $\Si$ are very special. In a way, they behave like round circles in $S^2_\infty(\BH^3)$.

\begin{lem} \label{rectangle} \cite[Lemma 3.2]{Co1} Let $R$ be a tall rectangle in $\Si$. Then there exists a unique minimal surface $\p$ in $\BHH$ with $\PI \p = \partial R$. 
\end{lem}

Furthermore, \cite{ST} gave an explicit description of the disk type minimal surface $\p$ \cite[Section 3]{Co1}.

The key component of our construction is the positive solution of following special case of Asymptotic Plateau Problem.

\begin{lem} \label{APP} \cite[Theorem 4.1]{Co1} [Tall Curves are Strongly Fillable] Let $\Gamma$ be a finite collection of disjoint, smooth Jordan curves in $\Si$ with $h(\Gamma)\neq \pi$. Then, there exists a complete, embedded, area minimizing surface $\Sigma$ in $\BHH$ with $\PI \Sigma = \Gamma$ if and only if $\Gamma$ is a tall curve.	
\end{lem}

Next lemma is an asymptotic regularity result for complete, embedded, area minimizing surfaces in $\BHH$.

\begin{lem} \label{surface} \cite[Lemma 7.6]{Co1} Let $\Sigma$ be a complete area minimizing surface in $\BHH$. Let $\overline{\Sigma}$ be the closure of $\Sigma$ in $\overline{\BHH}$, and let $\Gamma=\PI\Sigma$. If $\Gamma$ is a tall curve, then $\overline{\Sigma}$ is a surface with boundary.
\end{lem}

\begin{rmk} \label{regularityrmk} Notice that in the lemma above, everything is in $\C^0$ category. In \cite[Section 3]{KM}, Kloeckner and Mazzeo proved a stronger asymptotic regularity result for complete, embedded, minimal surfaces in $\BHH$ bounding  $\C^{k,\alpha}$ smooth curves in $\Si$.
\end{rmk}

The following classical result of geometric measure theory will be very useful for our construction.

\begin{lem} \label{AMSexist} \cite[Theorem 5.1.6 and 5.4.7]{Fe} [Existence and Regularity of Area Minimizing Surfaces] Let $M$ be a homogeneously regular, closed (or mean convex) $3$-manifold. Let $\gamma$ be a nullhomologous smooth curve in $\gamma$. Then, $\gamma$ bounds an area minimizing surface $\Sigma$ in $M$. Furthermore, any such area minimizing surface is smoothly embedded.	
\end{lem}

Now, we state the convergence theorem for area minimizing surfaces, which will be used throughout the paper. Note that we use convergence in the sense of Geometric Measure Theory, i.e. the convergence of rectifiable currents in the flat metric.

\begin{lem} \label{convergence} [Convergence] Let $\{\Sigma_i\}$ be a sequence of complete area minimizing surfaces in $\BHH$ where $\Gamma_i= \PI \Sigma_i$ is a finite collection of closed curves in $\Si$. If $\Gamma_i$ converges to a finite collection of closed curves $\wh{\Gamma}$ in $\Si$, then there exists a subsequence $\{\Sigma_{n_j}\}$ such that $\Sigma_{n_j}$ converges to an area minimizing surface $\wh{\Sigma}$ (possibly empty) with $\PI\wh{\Sigma}\subset\wh{\Gamma}$. In particular, the convergence is smooth on compact subsets of $\BHH$.
\end{lem}

\begin{pf}  Let $\Delta_n = \mathbf{B}_n(0)\times [-C,C]$ be convex domains in $\BHH$ where $\mathbf{B}_n(0)$ is the closed disk of radius $n$ in $\BH^2$ with center $0$, and $\wh{\Gamma}\subset \si\times (-C,C)$. For $n$ sufficiently large, consider the surfaces $S_i^n=\Sigma_i\cap \Delta_n$. Since the area of the surfaces $\{S^n_i\subset \Delta_n\}$ is uniformly bounded by $|\partial \Delta_n|$, and $\partial S^n_i$ can be bounded by using standard techniques. Hence, if $\{S^n_i\}$ is an infinite sequence, then we get a convergent subsequence of $\{S^n_i\}$ in $\Delta_n$ with {\em nonempty limit} $S^n$. $S^n$ is an area minimizing surface in $\Delta_n$ by the compactness theorem  for rectifiable currents (codimension-1) with the flat metric of Geometric Measure Theory (Lemma \ref{AMSexist}). By the regularity theory, the limit $S^n$ is a smoothly embedded area minimizing surface in $\Delta^n$.
	
If the sequence $\{S^n_i\}$ is an infinite sequence for infinitely many $n$, we get an infinite sequence of compact area minimizing surfaces $\{S^n\}$. Then, by using the diagonal sequence argument, we can find a subsequence of $\{\Sigma_i\}$ converging to an area minimizing surface $\wh{\Sigma}$ with $\PI \wh{\Sigma}\subset \wh{\Gamma}$ as $\Gamma_i\to \wh{\Gamma}$. Note also that for fixed $n$, the curvatures of $\{S^n_i\}$ are uniformly bounded by curvature estimates for area minimizing surfaces. Hence, with the uniform area bound, we get smooth convergence on compact subsets of $\BHH$. For further details, see \cite[Theorem 3.3]{MW}.
\end{pf}

\begin{rmk} Notice that in the lemma above, we can allow $\Gamma_i\subset \Si$ to be a collection of closed curves which may not be simple. Let $\Sigma_i$ be an area minimizing surface in $\BHH$ with $\PI\Sigma_i=\Gamma_i$. As $\Sigma_i$ is an area minimizing surface in $\BHH$, it must be embedded by the regularity of area minimizing surfaces. Hence, in such case, $\overline{\Sigma}_i$ may not be an embedded surface with boundary in $\overline{\BHH}$, even though $\Sigma_i$ is an open embedded surface in $\BHH$. Similarly, the limit $\wh{\Sigma}$ (if nonempty) is an open embedded surface in $\BHH$, even if $\PI\wh{\Sigma}\subset\wh{\Gamma}$ is not embedded. For the case $\Gamma\subset\Si$ is tall and embedded, see also Lemma \ref{surface}.
\end{rmk}

\section{Vertical Bridges at Infinity} \label{secbridge}

In this section, we  prove a bridge principle at infinity for sufficiently long vertical bridges. Then, by using these bridges, we  construct area minimizing surfaces of arbitrary topology in $\BHH$ in the next section.

\begin{defn} Let $\Gamma$ be  a collection of disjoint Jordan curves in $\Si$. If $\Gamma$ bounds a unique area minimizing surface $\Sigma$ in $\BHH$, we call $\Sigma$ {\em a uniquely minimizing surface}, and we call $\Gamma$ {\em a uni-curve}.
\end{defn}

\noindent {\em Notation and Setup:} Let $L_{\theta_0}$ be a vertical line in $\Si$, i.e. $L_{\theta_0}=\{\theta_0\}\times\BR$. Let $K_0>\pi$ be as in Lemma \ref{circles}.
Let $\Gamma$ be a smooth tall uni-curve in $\Si$ with $\Gamma\cap L_{\theta_0}=\emptyset$ and $h(\Gamma)>K_0$. Let $\Omega^\pm$ be the tall regions in $\Si$ with $\Gamma^c=\Omega^+\cup\Omega^-$ and $\partial \overline{\Omega^\pm}=\Gamma$.

Let $\alpha=\{\theta_1\}\times [c_1,c_2]$ be a vertical line segment in $\Si$ such that $\alpha\cap \Gamma = \partial \alpha$ and $\alpha \perp \Gamma$. Notice that $h(\Gamma)>K_0$   implies $c_2-c_1>K_0$, and $\alpha\cap \Gamma = \partial \alpha$ implies $\alpha\subset \overline{\Omega^+}$ or $\alpha\subset \overline{\Omega^-}$.

\begin{figure}[h]
	\begin{center}
		$\begin{array}{c@{\hspace{.5in}}c}

		\relabelbox  {\epsfysize=2in \epsfbox{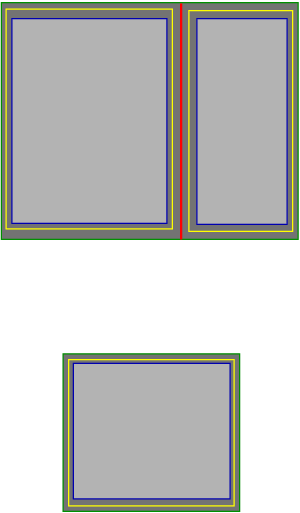}}   \endrelabelbox &
		
		\relabelbox  {\epsfysize=2in \epsfbox{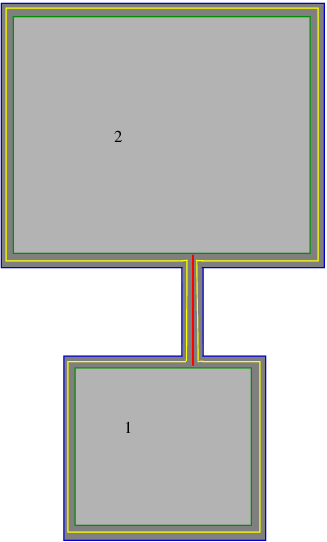}} \relabel{1}{ } \relabel{2}{ }  \endrelabelbox \\
		\end{array}$
		
	\end{center}
	\caption{\label{bridgefig} \footnotesize In the figure, $\Gamma=\partial \Omega^\pm$ is the green curves with two components. Light shaded regions (in the right) represent $\Omega^+$. In the left, we picture the case when the bridge $\alpha$ (red vertical line segment) is in $\Omega^+$. In the right, we picture the case when $\alpha$ is in $\Omega^-$. The family $\{\Gamma_t\}$ (yellow curves) foliate $\wh{N}$ (dark shaded region). Here, $\Gamma_\e \subset \partial \wh{N}$ is the blue curves.}
\end{figure}

Consider a small open neighborhood $N(\Gamma\cup\alpha)$ of $\Gamma\cup\alpha$ in $\Si$. If $\alpha\subset \Omega^+$, let $\wh{N}=N(\Gamma\cup\alpha)\cap \Omega^+$. If $\alpha\subset \Omega^-$, let $\wh{N}=N(\Gamma\cup\alpha) \cap \Omega^-$. In other words, we only take one side $\wh{N}$ of the open neighborhood $N(\Gamma\cup\alpha)$.  Foliate $\wh{N}$ by the smooth curves $\{\Gamma_t\ | \ t \in (0,\epsilon)\}$ with $\Gamma_\epsilon\subset \partial \wh{N}$, and $\Gamma_0 = \Gamma\cup \alpha$ (See Figure \ref{bridgefig}).  By taking a smaller neighborhood $N(\Gamma\cup\alpha)$ to start if necessary, we can assume that $\Gamma_t$ is a smooth tall curve for any $t$.

Let $S_\alpha$ be a thin strip along $\alpha$ in $\Si$. In particular, if $N(\alpha)$ is a small neighborhood of $\alpha$ in $\Si$, then $S_\alpha$ is the component of $N(\alpha)-\Gamma$ containing $\alpha$, i.e. $S_\alpha\sim [\theta_1-\delta,\theta_1+\delta]\times[c_1,c_2]$. In Figure \ref{bridgefig}, a tall curve $\Gamma$ with two components is pictured. In the left figure, the bridge $\alpha$ is in $\Omega^+$, while in the right, $\alpha$ is in $\Omega^-$. Notice that if $\partial \alpha$ is in the same component of $\Gamma$, then $\sharp(\Gamma_t)=\sharp(\Gamma)+1$ where $\sharp(.)$ represents the number of components (Figure \ref{bridgefig} Left). Similarly, if $\partial \alpha$ is in the different components of $\Gamma$, then $\sharp(\Gamma_t)=\sharp(\Gamma)-1$ (Figure \ref{bridgefig} Right).

Now, consider the upper half plane model for $\BH^2\simeq \{(x,y)\ | \ y>0\}$. Without loss of generality, let $\theta_0\in S^1_\infty(\BH^2)$ correspond to the point at infinity in the upper half plane model. We  use the upper half space model for $\BHH$ with the identification $\BHH=\{(x,y,z) \ | \ y>0\}$ where $\BH^2$ corresponds the $xy$-half plane, and $\BR$ corresponds to $z$ coordinate. Hence, the $xz$-plane will correspond to $\Si$. By using the isometries of the hyperbolic plane and the translation along $\BR$ direction, we  assume that $\theta_1\in S^1_\infty(\BH^2)$ will correspond to $0$, and the vertical line segment $\alpha\subset \Si$ above will have $\alpha=\{(0,0)\}\times [c_1,c_2]$ and $S_\alpha\sim [-\delta,\delta]\times\{0\}\times [c_1,c_2]$ in $(x,y,z)$ coordinates.

With this notation, we can state the bridge principle at infinity for vertical bridges in $\Si$ as follows.

\begin{thm} [Vertical Bridges at Infinity] \label{bridge} Let $\Gamma$ be a tall uni-curve with $h(\Gamma)>K_0$ as above. Define $\alpha, \Gamma_t, S_\alpha$ accordingly as described above. Let $\Sigma$ be the uniquely minimizing surface in $\BHH$ where $\PI \Sigma=\Gamma$. Assume also that $\overline{\Sigma}$ has finite genus. Then, there exists a sufficiently small $t>0$ such that $\Gamma_t$ bounds a unique area minimizing surface $\Sigma_t$ where $\Sigma_{t}$ is homeomorphic to $\Sigma\cup S_\alpha$, i.e. $\Sigma_t \simeq \Sigma\cup S_\alpha$.
\end{thm}

\noindent {\em Outline of the Proof:} Let $\Gamma_t\to(\Gamma\cup\alpha)$ as above. Let $\Sigma_t$ be the area minimizing surface in $\BHH$ with $\PI\Sigma_t=\Gamma_t$.  Intuitively, for sufficiently large $n>0$, we want to show that $\Sigma_{t_n}$ is just $\Sigma$ with a thin strip along $\alpha$, where thin strip vanishes as $n\to \infty$. We split the proof into 4 steps. 
In Step 1, we blow up the sequence $\{\Sigma_{t_n}\}$ and show that the limit $T=\lim \Sigma_{t_n}$ cannot contain the vertical segment $\alpha$. In Step 2, we  show that $\Sigma_t$ does not develop any genus near the asymptotic boundary. In Step 3, we show that $\Sigma_t \simeq \Sigma\cup S_\alpha$ for $t$ sufficiently close to $0$. Finally in Step 4, by using generic uniqueness, we show that we can choose $t>0$ such that $\Gamma_t$ bounds a unique area minimizing surface $\Sigma_t$.

\begin{pf} First, by Lemma \ref{APP}, for any $\Gamma_t\subset \Si$, there exists an area minimizing surface $\Sigma_t$ with $\PI \Sigma_t = \Gamma_t$.

As $t_n\searrow 0$, $\Gamma_{t_n}\to\Gamma\cup\alpha$. Since $\Gamma_{t_n}$ is a tall curve, there exists an area minimizing surface $\Sigma_{t_n}$ in $\BHH$ with $\PI\Sigma_{t_n}=\Gamma_{t_n}$ by Lemma \ref{APP}. By Lemma \ref{convergence}, there exists a convergent subsequence, say $\Sigma_n$, converging to an area minimizing surface $T$ with $\PI T\subset\Gamma\cup\alpha$. Since $\Gamma\cup\alpha$ is a tall curve, the limit $T$ is nonempty by the proof of Lemma \ref{APP}.

Now, we claim that $\PI T = \Gamma$. In other words, the limit area minimizing surface $T$ with $\PI T\subset\Gamma\cup\alpha$ cannot have the vertical segment $\alpha$ in its asymptotic boundary. Then, since $\Gamma$ bounds a unique area minimizing surface $\Sigma$, $\PI T = \Gamma$ would imply $T=\Sigma$.

\vspace{.2cm}

\noindent {\bf Step 1:} $\PI T = \Gamma$.

\vspace{.2cm}

{\em Proof of Step 1:} By above, we know that $\PI T \subset \Gamma\cup \alpha$. Note that by Lemma \ref{surface}, $\overline{T}=T\cup\PI T$ is a surface with boundary in $\overline{\BHH}$. As $\overline{\BHH}$ is topologically a closed ball, $\overline{T}$ is separating in $\overline{\BHH}$.

Assume that there is a point $p\in \alpha-\partial\alpha$ such that $p\in \PI T$. By using the notation and the upper half space model described before the theorem, recall that $\alpha=\{(0,0)\}\times [c_1,c_2]$, and without loss of generality, assume $p=(0,0,0)\in \alpha\subset \Si$. Consider the hyperbolic plane $P=\BH^2\times \{0\}=\{(x,y,0) \ | \ y>0\}$ in $\BHH$.  Let $\gamma_i$ be the geodesic arc in $P$ with $\PI \gamma_i =\{(-r_i,0,0),(+r_i,0,0)\}$ where $r_i\searrow 0$. Let $\U_i=\gamma_i\times[-\e_0,+\e_0]$ for some fixed $\e_0$. Then, since $p\in \PI T$, $T\cap \U_i\neq \emptyset$ for $i>N_0$ for some $N_0$. Let $q_i\in T\cap \U_i$ for any $i>N_0$. 

Now, let $\varphi_i$ be the isometry of $\BHH$ with $\varphi_i(x,y,z)= (\frac{1}{r_i}x,\frac{1}{r_i}y, z)$. Define a sequence of area minimizing surfaces $T_i=\varphi_i(T)$. Let $\wh{\gamma}$ be the geodesic in $P$ with $\PI \wh{\gamma}=\{(-1,0,0),(1,0,0)\}$. Hence, by construction, $\varphi_i(\gamma_i)=\wh{\gamma}$ and $\varphi_i(\U_i)=\wh{\U}=\wh{\gamma}\times [-\e_0,+\e_0]$ for any $i>0$. Let $\wh{q}_i=\varphi_i(q_i)$ for any $i>N_0$. Then, $\wh{q}_i\subset T_i\cap\wh{\U}$ for any $i>N_0$. 
Again by using Lemma \ref{convergence}, we get a subsequence of $\{T_i\}$ which converges to an area minimizing surface $\wh{T}$. Let $R^+$ and $R^-$ be two tall rectangles in opposite sides of $\alpha$ disjoint from $\Gamma\cup\alpha$, and let $P^\pm$ be the unique area minimizing surfaces with $\PI P^\pm=\partial R^\pm$. By Lemma \ref{disjoint} and Remark \ref{disjointrem}, $T_i\cap P^\pm=\emptyset$ for any $i$. Let  $\eta$ be the finite segment in $\wh{\gamma}$ with $\partial \eta\subset P^+\cup P^-$. Let $\wh{\V}=\eta\times [-\e_0,+\e_0]$. Then, $\{\wh{q}_i\}\subset \wh{\V}\subset\wh{\U}$. As $\wh{\V}$ is compact, $\{\wh{q}_i\}$ has a convergent subsequence. This implies $\wh{T}\cap \wh{\V}\neq \emptyset$.  This proves that the limit area minimizing surface $\wh{T}$ does not escape to infinity. Furthermore, in above construction, we can choose $\wh{\U}$ as close as we want to infinity $\{(0,0)\}\times[-\e_0,\e_0]$, and we can choose $\e_0>0$ as small as we want, we conclude that $p\in \PI \wh{T}$, too.

Now, by the construction of the sequence $\{T_i\}$, $\wh{T}$ and hence $\PI\wh{T}$ are invariant by the isometry $\varphi_\lambda(x,y,z)= (\lambda x,\lambda y, z)$. Notice that the isometry $\varphi_\lambda$ fixes only the points $(0,0)$ and $\infty$ in $\si$ and the horizontal lines $L_i=\{(t,0,c_i) \ | \ t\in\BR\}$ in $\Si$.
This implies $\PI\wh{T}\subseteq \wh{\Gamma}$ where $\wh{\Gamma}\subset\alpha\cup L_1\cup .. L_{m_1}\cup \beta_1\cup .. \beta_{m_2}$ where $\beta_j$ is a vertical line segments with $x$-coordinate $0$. In particular, in the cylindrical model for $\BHH$, $\wh{\Gamma}\subset \alpha\bigcup_{i=1}^{m_1} \gamma_{c_i}\bigcup_{j=1}^{m_2} \beta_j\bigcup_{k=1}^{m_3} \wh{\beta}_k$ where $\gamma_{c_i}=\si\times\{c_i\}$ is the horizontal circle corresponding to $L_i$ in $\Si$. $\beta_j=\{\theta_1\}\times[c^-_j,c^+_j]$ and $\wh{\beta}_k=\{\theta_0\}\times[d^-_k,d^+_k]$ where $\theta_0\sim \infty$ and $\theta_1\sim (0,0)$ in the upper half space model. Since $h(\Gamma)>K_0$, then $c^+_j-c^-_j>K_0>\pi$ and $d^+_k-d^-_k>K_0>\pi$ by construction. This implies the area minimizing surface $\wh{T}$ satisfies the conditions of Lemma \ref{circles} in the appendix. By the lemma, we conclude that $\PI\wh{T}\subset \bigcup_{i=1}^{m_1} \gamma_{c_i}$, i.e. $\PI\wh{T}$ is a collection of horizontal circles in $\Si$, and cannot have any vertical line segments like $\alpha$. However, this gives a contradiction as $p\in \PI\wh{T}$.  Step 1 follows. \hfill $\Box$

Now, we show that $\Sigma_t$ does not develop genus near the asymptotic boundary.

\vspace{.2cm}

\noindent {\bf Step 2:}  There exists $a_\Gamma>0$ such that for sufficiently large $n$, $\Sigma_n\cap \R_{a_\Gamma}$ has no genus, i.e. $\Sigma_n\cap \R_{a_\Gamma} \simeq \Gamma_n\times (0,a_\Gamma)$.

\vspace{.2cm}

\noindent {\em Proof of the Step 2:} Assume on the contrary that for  $a_n\searrow0$, there exists a subsequence $\Sigma_n\cap \R_{a_n}$ with positive genus. Recall that by Lemma \ref{surface}, $\overline{\Sigma}_n=\Sigma_n\cup\Gamma_n$ is a surface with boundary in $\overline{\BHH}$ and separating in $\overline{\BHH}$. Let $\Delta_n$ be the component of $\overline{\BHH}-\overline\Sigma_n$ which contains the bridge $\alpha$. Since $\Sigma_n\cap \R_{a_n}$ has positive genus, then $\Delta_n\cap \R_{a_n}$ must be a nontrivial handlebody, i.e. not a $3$-ball. Hence, there must be a point $p_n$ in $\Sigma_n\cap \R_{a_n}$ where the normal vector $v_{p_n}=\langle 0,1,0\rangle$ pointing {\em inside} $\Delta_n$ by Morse Theory. By genericity of Morse functions, we can modify the $\infty$ point in $\PI\BH^2$ if necessary, to get $y$-coordinate as a Morse function.


Let $p_n=(x_n,y_n,z_n)$. By construction $y_n\to 0$ as $y_n<a_n$. Consider the isometry $\psi_n(x,y,z)=(\frac{x-x_n}{y_n},\frac{y}{y_n},z-z_n)$ which is a translation by $-(x_n,0,0)$ first by a parabolic isometry of $\BH^2$, and translation by $-(0,0,z_n)$ in $\BR$ direction. Then, by composing with the hyperbolic isometry $(x,y,z)\to (\frac{x}{y_n},\frac{y}{y_n},z)$, we get the isometry $\psi_n$ of $\BHH$. Then, consider the sequence of area minimizing surface $\Sigma_n'=\psi_n(\Sigma_n)$ and $p_n'=\psi_n(p_n)=(0,1,0)$. Let $\Gamma_n'=\psi_n(\Gamma_n)=\PI \Sigma_n'$. After passing to a subsequence, we get the limits $\Sigma_n'\to \Sigma'$, $p_n'\to p'=(0,1,0)\in \Sigma'$, and $\Gamma_n'\to\Gamma'$. Note also that by construction the normal vector to area minimizing surface $\Sigma'$ at $p'$ is $v_{p_n}\to v_p'=<0,1,0>$ pointing inside $\Delta'$.

Consider $\Gamma'=\lim \Gamma_n'$. Let $l_z$ be the $z$-axis in $\Si$, i.e. $l_z=\{(0,0,t) \ | \ t\in\BR\}$. Let $\Gamma'\cap l_z=\{ (0,0,c_1), (0,0,c_2), ..,(0,0,c_k)\}$. Notice that as $h(\Gamma)>K_0$ , $|c_i-c_j|>K_0$ for any $i\neq j$. Recall that $\partial \alpha=\{(0,0,c_1),(0,0,c_2)\}$. Note that by Lemma \ref{circles}, $\Gamma'$ cannot have a vertical line segment $\alpha_j=\{(0,0)\}\times[c_j^-,c_j^+]$. Hence, by construction of $\Gamma_n'$, we get $\Gamma'= \beta\cup L_{c_3}\cup ..\cup L_{c_k}$ where $L_{c_i}$ is the horizontal line in $\Si$ with $L_{c_i}$, and $\beta$ is the component of $\Gamma'$ near $\alpha$ (See Figure \ref{gamma'} left). In particular, in cylinder model for $\BHH$, $L_{c_i}$ is the horizontal circle $\gamma_{c_i}=\si\times\{c_i\}$ in $\Si$, and $\beta$ is a tall rectangle $\beta=\partial R$ where $R=[\delta,2\pi-\delta]\times [c_1,c_2]$ assuming $\alpha=\{0\}\times[c_1',c_2']$ (See Figure \ref{gamma'} right). Note that as $\psi_n$ is only translating in $z$-direction, $c_1-c_2=c_1'-c_2'$. Here,  the limit area minimizing surface $\Sigma'$ is nonempty, as $(0,0,1)\in \Sigma'$ by construction. $\delta$ depends on the comparison between $y_n\searrow 0$ and $d(\Gamma_n, \alpha)\searrow 0$. As $\Sigma'$ does not escape infinity, we make sure that such a $\delta<\pi$ exists. Indeed, $\delta>0$ can be explicitly computed by using the fact that there is a unique minimal surface $P_\beta$ in $\BHH$ containing (0,0,1) with $\PI P_\beta = \beta=\partial R$ by Lemma \ref{rectangle} as $R$ is a tall rectangle.

\begin{figure}[t]
	\begin{center}
		$\begin{array}{c@{\hspace{.5in}}c}

		\relabelbox  {\epsfxsize=2in \epsfbox{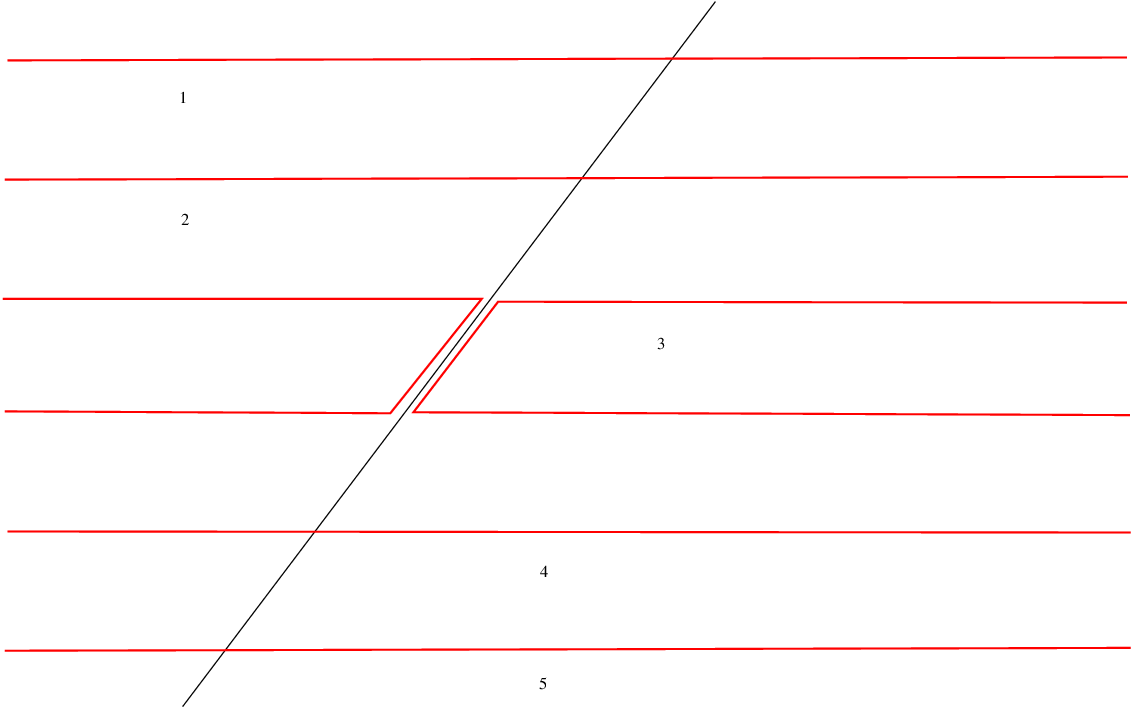}}  \relabel{1}{\tiny $L_{c_6}$} \relabel{2}{\tiny  $L_{c_5}$} \relabel{3}{\tiny  $\beta$} \relabel{4}{\tiny  $L_{c_3}$} \relabel{5}{\tiny  $L_{c_4}$}   \endrelabelbox &
		
		\relabelbox  {\epsfysize=1.7in \epsfbox{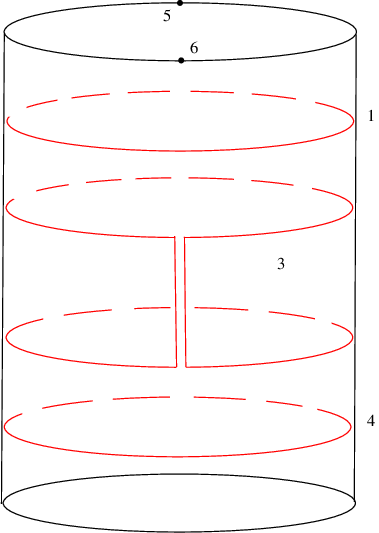}} \relabel{1}{\tiny  $\gamma_{c_5}$} \relabel{3}{\tiny  $\beta$} \relabel{4}{\tiny  $\gamma_{c_3}$} \relabel{5}{\tiny $\pi$}  \relabel{6}{\tiny $0$}  \endrelabelbox \\
	\end{array}$
	
\end{center}
\caption{\label{gamma'} \footnotesize $\Gamma'\subset \Si$ is pictured in upper half space model and cylinder model for $\BHH$. }
\end{figure}

$\Sigma'$ bounds a unique area minimizing surface with $\Sigma'= P_\beta\cup P_{c_1} \cup ..\cup P_{c_k}$ where $P_\beta$ is the unique area minimizing surface with $\PI P_\beta= \beta$ by Lemma \ref{rectangle}, and $P_{c_i}$ is the horizontal plane $\BH^2\times \{c_i\}$ in $\BHH$ with $\PI P_{c_i}= \gamma_{c_i}$. This is because $|c_i-c_j|>\pi$, there is no connected minimal surface with asymptotic boundary contains more than one component of $\Gamma'$. In particular, if there was a connected area minimizing surface $Y$ with $\PI Y\supset \gamma_{c_1}\cup\gamma_{c_2}$ with $c_2-c_1>\pi$, one can place a minimal catenoid $\C$ with $\PI \C=\gamma_{c_1'}\cup\gamma_{c_2'}$ where $c_1'>c_1$ and $c_2'<c_2$ with $c_2'-c_1'=\pi-\e$ so that $\C\cap S=\emptyset$ (\cite[Section 7.1]{Co1}). Then by using an hyperbolic isometry $\varphi_t$, one can push $\C$ towards $S$ horizontally. As $S$ is connected, there must be a first point of touch, which gives a contradiction with maximum principle. This implies each component of $\Gamma'$ bounds a component of $\Sigma'$. Since each component is uniquely minimizing, $\Sigma'$ is a uniquely minimizing surface with $\PI\Sigma'=\Gamma'$.

Hence, by construction $p'=(0,0,1)$ is on $P_\beta$ component of $\Sigma'$. Recall that the normal vector $v_p'=\langle 0,1,0\rangle$ points inside of $\Delta'$ which is the component of $\BHH-\Sigma'$ containing $\alpha$. However, $P_\beta$ is a plane, and the normal vector $v_p'$ points outside of $\Delta'$ not inside. This is a contradiction. Step 2 follows. \hfill $\Box$

\vspace{.2cm}

\noindent {\bf Step 3:} For sufficiently small $t>0$, $\Sigma_{t}$ is homeomorphic to $\Sigma\cup S_\alpha$.

\vspace{.2cm}

\noindent {\em Proof of Step 3:} Assume that for $\e_n\searrow 0$, there exists $0<t_n<\e_n$ such that $\Sigma_{t_n}$, say $\Sigma_n$, is not homeomorphic to $\wh{\Sigma}=\Sigma\cup S_\alpha$. Since the number of ends are same, this means $\Sigma_n$ and $\wh{\Sigma}$ have different genus.


Let $\R_a=\{0\leq y \leq a\}$ in $\overline{\BHH}$ be as in Step 2. Let $\mathcal{K}_a=\{y\geq a\}$ and let $\Sigma^a=\Sigma\cap \mathcal{K}_a$. Then, since $\Sigma_n\to \Sigma$ converge smoothly on compact sets,  $\Sigma^a_n\to \Sigma^a$ smoothly. Hence, by Gauss-Bonnet, $\Sigma^a_n$ and $\Sigma^a$ must have same genus. By Step 2, this implies for sufficiently large $n$, $\Sigma_n$ and $\Sigma$ must have the same genus. However, this contradicts with our assumption that $\Sigma_n$ and $\Sigma$ have different genus for any $n$. Therefore, this implies that for sufficiently small $\epsilon'>0$, $\Sigma_t$ is homeomorphic to $\Sigma\cup S_\alpha$ for $0<t<\epsilon'$. Step 3 follows. \hfill $\Box$

\vspace{.2cm}

\noindent {\bf Step 4:} For all but countably many $0<t<\epsilon'$, $\Gamma_t$ bounds a unique area minimizing surface in $\BHH$.

\vspace{.2cm}

\noindent {\em Proof of Step 4:}  We  adapt the proof of Theorem \ref{uniq} to this case. The family of tall curves $\{\Gamma_t\ | \ t \in (0,\epsilon)\}$ foliates $\wh{N}$ where $\partial \wh{N}=\Gamma_\epsilon\cup \Gamma$, and $\Gamma_0 = \Gamma\cup \alpha$. In particular, for any $0<t_1<t_2<\e$, $\Gamma_{t_1}\cap \Gamma_{t_2}=\emptyset$. If $\Sigma_t$ is an area minimizing surface in $\BHH$, then $\Sigma_{t_1}\cap\Sigma_{t_2}=\emptyset$ too, by Lemma \ref{disjoint}. By Lemma \ref{canonical}, if $\Gamma_s$ does not bound a unique area minimizing surface $\Sigma_s$, then we can define two disjoint canonical minimizing  $\Sigma^+_s$ and $\Sigma^-_s$ with $\PI \Sigma^\pm_s = \Gamma_s$. Then, by the proof of Theorem \ref{uniq}, for all but countably many $s\in[0,\epsilon']$, $\Gamma_s$ bounds a unique area minimizing surface. Step 4 follows. \hfill $\Box$




Steps 1-3 implies the existence of $\e'>0$ such that any $\Sigma_t$ with $\PI\Sigma_t=\Gamma_t$ for $t\in (0,\e')$ is homeomorphic to $\Sigma\cup S_\alpha$. Step 4 implies the generic uniqueness for the family $\{\Gamma_t\mid t\in(0,\e')\}$. Hence, Step 1-4 together implies the existence of smooth curve $\Gamma_t$ with $t\in(0,\e')$, where $\Gamma_t$ bounds a {\em unique} area  minimizing surface $\Sigma_t$, and $\Sigma_t$ has the desired topology, i.e. $\Sigma_t\simeq \Sigma\cup S_\alpha $. The proof of the theorem follows.
\end{pf}

\section{Minimal Surfaces of Arbitrary Topology in $\BHH$} \label{secsurfaceconstruction}

In this section, we prove any open orientable surface can be embedded in $\BHH$ as an area minimizing surface. First, we show a simple construction for finite topology case. Then, we finish the proof by giving a very general construction for infinite topology case.

\subsection{Surfaces with Finite Topology} \label{chi} \

\vspace{.2cm}

While our main result later implies both finite or infinite topology orientable surface, we start with a very simple construction for surfaces with finite topology as a warm-up. In particular, by using vertical bridges as $1$-handles, we give a construction of an area minimizing surface $\Sigma_k^g$ of genus $g$ with $k$ ends.

\vspace{.2cm}

\noindent {\em Euler Characteristics:} Recall that if $T^g_k$ is an orientable surface of genus $g$, and $k$ boundary components, then $\chi(T^g_k)=2-2g-k$. Adding a bridge (a 1-handle in topological terms) to a surface decreases the Euler Characteristics by one. On the other hand, if you add a bridge to a surface where the endpoints of the bridge are in the same boundary component, then the number of boundary components increases by one. If you add a bridge whose endpoints are in the different boundary components, then the number of boundary components decreases by one (See Figure \ref{bridgefig}).

Now, adding a bridge to the same boundary component of a surface would increase the number of ends. In other words, let $S_{n+1}$ obtained from $S_n$ by attaching a bridge ($1$-handle) to $S_n$ whose endpoints are in the same component of $\partial S_n$. Then, $\chi(S_{n+1})=\chi(S_n)-1$, $g(S_n)=g(S_{n+1})$ and $\sharp (\partial S_{n+1})=\sharp (\partial S_n)+1$ where $\sharp$ is the number of components.

If we want to increase the genus, first add a bridge $\alpha_n$ whose endpoints are in the same component of $\partial S_n$, and get $S_n'\simeq S_n\natural S_{\alpha_n}$ where $S_n\natural S_{\alpha_n}$ represents the surface obtained by adding a bridge (thin strip) to $S_n$ along $\alpha_n$. Then, by adding another bridge $\alpha'_n$ whose endpoints are in different components of $S_n'$, one get $S_{n+1} \simeq S_n'\natural S_{\alpha'_n}$.  Hence, $\chi(S_{n+1})=\chi(S_n)-2$, and the number of boundary components are same. This implies if $S_n \simeq T^g_k$, then $S_{n+1}\simeq T^{g+1}_k$. This shows that $S_{n+1}$ is obtained by attaching a cylinder with handle to $S_n$, i.e. $S_{n+1}-S_n$ is a cylinder with handle.


\begin{figure}[t]
	\begin{center}
		$\begin{array}{c@{\hspace{.4in}}c}

		\relabelbox  {\epsfysize=2in \epsfbox{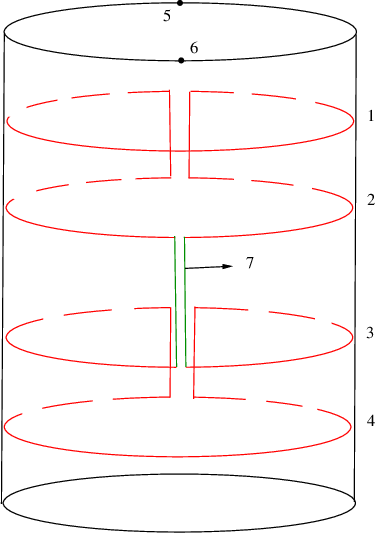}} \relabel{1}{\tiny $2K_0$} \relabel{2}{\tiny $K_0$} \relabel{3}{\tiny $-K_0$} \relabel{4}{\tiny $-2K_0$} \relabel{5}{\tiny $\pi$}  \relabel{6}{\tiny $0$} \relabel{7}{\small $\tau$}  \endrelabelbox &
		
		\relabelbox  {\epsfysize=2in \epsfbox{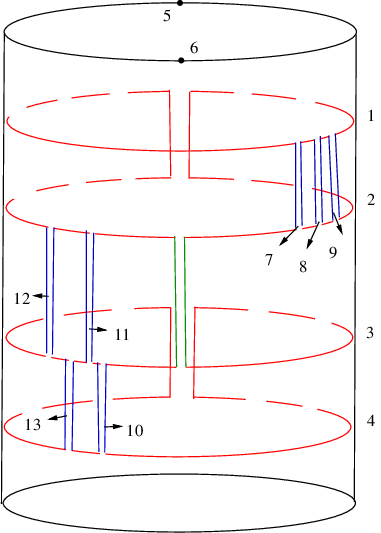}} \relabel{1}{\tiny $2K_0$} \relabel{2}{\tiny $K_0$} \relabel{3}{\tiny $-K_0$} \relabel{4}{\tiny $-2K_0$} \relabel{5}{\tiny $\pi$}  \relabel{6}{\tiny $0$} \relabel{7}{\tiny $\beta_1$} \relabel{8}{\tiny $\beta_2$} \relabel{9}{\tiny $\beta_3$} \relabel{10}{\tiny $\zeta_1$} \relabel{11}{\tiny $\zeta_1'$} \relabel{12}{\tiny $\zeta_2'$} \relabel{13}{\tiny $\zeta_2$} \endrelabelbox \\
	\end{array}$
	
\end{center}
\caption{\label{sigma} \footnotesize In the left, we have the tall curve $\Gamma_1$ which bounds the area minimizing surface $\Sigma_1\sim P^+\cup P^-\natural S_\tau$. In the right, we first add bridges $\beta_1, .., \beta_k$ to $\Sigma$ to increase the number of ends by $k$ (here for  $k=3$). Then, we add $g$ pairs of bridges $\zeta_1,\zeta_1',...,\zeta_g,\zeta_g'$ to increase the genus (here $g=2$). Hence, $\Sigma$ is a genus 2 surface with 4 ends.}
\end{figure}

\vspace{.2cm}

\noindent {\em Construction for finite topology surfaces:} There is a very elementary construction for open orientable surfaces of finite topology as follows: Let $S$ be open orientable surface of genus $g$ and $k$ ends. Construct the area minimizing surface $\Sigma_1$ which is topologically a disk as in Figure \ref{sigma}-Left. For $k+1$ ends, add $k$ vertical bridges $\beta_1,\beta_2, ..,\beta_k$ to $\Sigma_1$ as in the Figure \ref{sigma}-Right. Then, for genus $g$, add $g$ pairs of vertical bridges $\zeta_i$ and $\zeta_i'$ successively as in Figure \ref{sigma}-Right. Hence, the final surface $\Sigma$ is an area minimizing surface of genus $g$ and $k+1$ ends. Furthermore, $\overline{\Sigma}$ is a compact embedded surface with boundary in $\overline{\BHH}$ by Lemma  \ref{surface}.

\subsection{Surfaces with Infinite Topology} \

\vspace{.2cm}

Now, we prove any open orientable surface (finite or infinite topology) can be embedded in $\BHH$ as an area minimizing surface. In this part, we  mainly follow the techniques in \cite{MW} and \cite{Co2}. In particular, for a given surface $S$, we  start with a compact exhaustion of $S$, $S_1\subset S_2\subset ...\subset S_n \subset ...$, and by using the bridge principle proved in the previous section, we inductively construct the area minimizing surface with the desired topology.

\begin{figure}[h]

\relabelbox  {\epsfxsize=3.5in

\centerline{\epsfbox{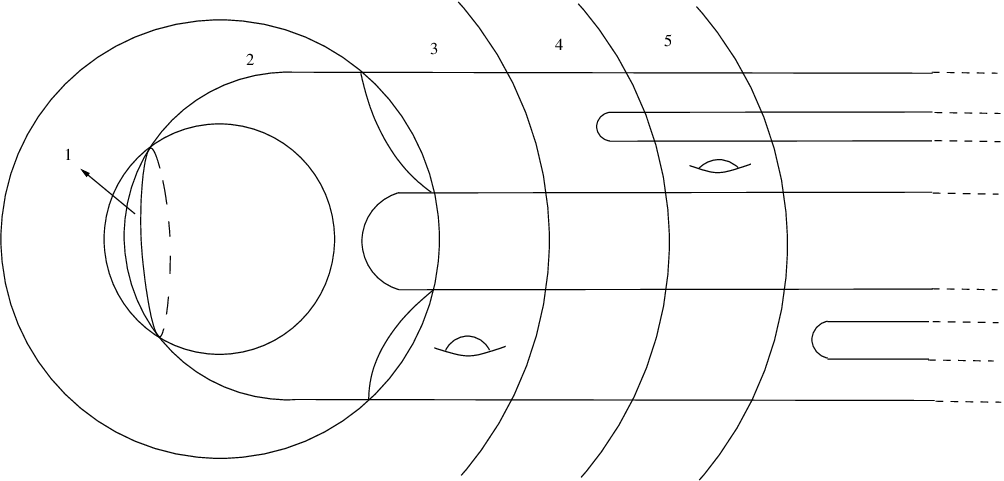}}}

\relabel{1}{$S_1$}

\relabel{2}{$S_2$}

\relabel{3}{$S_3$}

\relabel{4}{$S_4$}

\relabel{5}{$S_5$}

\endrelabelbox

\caption{\label{simpleexhaustion} \small In the simple exhaustion of $S$, $S_1$ is a disk, and $S_{n+1}-S_n$ contains a unique nonannular part, which is a pair of pants (e.g. $S_4-S_3$), or a cylinder with a handle (e.g. $S_3-S_2$). }

\end{figure}

In particular, by \cite{FMM}, for any open orientable surface $S$, there exists a simple exhaustion. A simple exhaustion $S_1\subset S_2\subset ...\subset S_n \subset ...$ is the compact exhaustion with the following properties: $S_1$ is a disk, and $S_{n+1}-S_n$ would contain a {\em unique nonannular piece} which is either a cylinder with a handle, or a pair of pants by \cite{FMM} (See Figure \ref{simpleexhaustion}).

First, we need a lemma which will be used in the construction.

\begin{lem} \label{disjrec} Let $R=[-1,1]\times[-4\pi, 4\pi]$ and $R_c=[-c,c]\times[-2\pi, 2\pi]$ be rectangles in $\Si$ where $0<c<1$. Let  $\gamma= \partial R$, $\gamma_c = \partial R_c$ and $\Gamma_c=\gamma\cup\gamma_c$. Then, there exists $\rho>0$ such that for any $0<c\leq\rho$, the area minimizing surface $\Sigma_c$ with $\PI \Sigma_c=\Gamma_c$ is $P\cup P_c$ where $P$ and $P_c$ are the unique area minimizing surfaces with $\PI P = \gamma$ and $\PI P_c = \gamma_c$.
\end{lem}

\begin{pf} If the area minimizing surface $\Sigma_c$ is not connected, then it must be $P\cup P_c$ because the rectangles $\gamma$ and $\gamma_c$ bounds a unique area minimizing surfaces $P$ and $P_c$ respectively by Lemma \ref{rectangle}. Hence, we assume on the contrary that the area minimizing surface $\Sigma_c$ is connected for any $0<c<1$. We abuse the notation and say $\Sigma_n= \Sigma_{\frac{1}{n}}$. Consider the sequence $\{\Sigma_n\}$. By Lemma \ref{convergence}, we get a convergent subsequence, and limiting area minimizing surface $\Sigma$ with $\PI \Sigma \subset \gamma\cup \beta$ where $\beta$ is the vertical line segment $\{0\}\times[-2\pi,2\pi]$.

Let $Q=[-\frac{1}{2},\frac{1}{2}]\times [-3\pi,3\pi]$ be another rectangle in $\Si$, and let $T$ be the unique area minimizing surface in $\BHH$ with $\PI T = \partial Q$. Since by assumption, $\Sigma_n$ is connected, and $\overline{T}$ separates the boundary components of $\Sigma_n$, $\gamma_n$ and $\gamma$, then $T\cap \Sigma_n\neq \emptyset$ for any $n>2$. By construction, this implies $\Sigma\cap T\neq \emptyset$. 

As $\PI \Sigma \subset \gamma\cup \beta$, we have two cases. Either $\PI\Sigma=\gamma$ or $\PI \Sigma = \gamma\cup \beta$. If $\PI\Sigma=\gamma$, then $\gamma$ bounds a unique area minimizing surface $P$. In other words, $\Sigma$ must be $P$ and $P\cap T =\emptyset$. This is a contradiction.

If $\PI \Sigma = \gamma\cup \beta$, we get a contradiction as follows. Let's go back to cylinder model of $\BHH$. Then, we can represent $\gamma=\partial R$ where $R= [-\theta_1,\theta_1]\times[-4\pi,4\pi]$ for some $\theta_1\in (0,\pi)$, and $\beta=\{0\}\times[-2\pi,2\pi]$ in $\Si$. Let $\varphi_t$ be the isometry of $\BHH$ corresponding to $\varphi_t(x,y,z)=(tx,ty,z)$ in upper half space model of $\BHH$. In particular, $\{\pi\}\times\BR$ represents the point at infinity, and $\varphi_t$ pushes every point in $\BHH$ from $\{0\}\times\BR$ to $\{\pi\}\times\BR$ in the Poincare disk model. Let $\Sigma^n=\varphi_n(\Sigma)$. Again, by Lemma \ref{convergence}, we get a limit area minimizing surface $\wh{\Sigma}$ where $\PI\wh{\Sigma}\subset \Gamma^+\cup\Gamma^-\cup\beta\cup \alpha$ where $\Gamma^\pm=\si\times\{\pm4\pi\}$ and $\alpha=\{\pi\}\times[-4\pi,4\pi]$. 

We claim that $\wh{\Sigma}$ is nonempty, and furthermore, $\PI\wh{\Sigma}= \Gamma^+\cup\Gamma^-\cup\beta\cup \alpha$. Since the original $\Sigma$ is connected by assumption, $\Sigma\cap \BH^2\times \{c\}$ contains an infinite curve $l_c$ with $\PI l_c=\{(0,c),(\theta_1,c)\}$ where $c\in(-2\pi,2\pi)$. Then, $\varphi_n(l_c)=l^n_c\subset\Sigma^n\cap \BH^2\times \{c\}$, and $l^n_c$  converges to a line $\wh{l}_c\subset \wh{\Sigma}\cap \BH^2\times \{c\}$ with $\PI\wh{l}_c=\{(0,c),(\pi,c)\}$. This shows $\PI\Sigma= \Gamma^+\cup\Gamma^-\cup\beta\cup \alpha$.

Finally, let $\C$ be the Daniel's parabolic catenoid with $\PI \C= \lambda^+\cup\lambda^-\cup \tau$ where $\lambda^+=\si\times \{\frac{7\pi}{2}\}$, $\lambda^-=\si\times \{\frac{5\pi}{2}\}$, and $\tau=\{\pi\}\times[\frac{5\pi}{2},\frac{7\pi}{2}]$. As $\PI\C$ is invariant by $\varphi_t$, $\C_t\varphi_t(\C)$ is also parabolic catenoid with $\PI \C_t=\PI\C$. Furthermore, for sufficiently small $\e>0$, $\C_\e$ is very close to asymptotic cylinder $\Si$. Hence, we can choose sufficiently small $\e>0$ with $\C_\e\cap \wh(\Sigma)=\emptyset$. Then, by pushing $\C_\e$ towards $\wh{\Sigma}$ via isometries $\varphi_t$, we  get a first point of touch $\C_{t_0}$ with $\wh{\Sigma}$ which contradicts to the maximum principle. The proof follows.
\end{pf}

Now, we are ready to prove the existence result for properly embedded area minimizing surfaces in $\BHH$ with arbitrary topology.

\begin{thm} \label{main} Any open orientable surface $S$ can be embedded in $\BHH$ as an area minimizing surface $\Sigma$.
\end{thm}

\begin{pf} Let $S$ be an open orientable surface. Now, we inductively construct an area minimizing surface $\Sigma$ in $\BHH$ which is homeomorphic to $S$. Let $S_1\subset S_2\subset ...\subset S_n \subset ...$ be a simple exhaustion of $S$, i.e. $S_{n+1}-S_n$ contains a unique nonannular piece which is either a cylinder with a handle, or a pair of pants.

By following the simple exhaustion, we  define a sequence of area minimizing surfaces $\Sigma_n$ so that $\Sigma_n$ is homeomorphic to $S_n$, i.e. $\Sigma_n\simeq S_n$. Furthermore, the sequence $\Sigma_n$  induces the same simple exhaustion for the limiting surface $\Sigma$. Hence, we get an area minimizing surface $\Sigma$ which is homeomorphic to the given surface $S$.

Now, we  follow the idea described in Remark \ref{chi}. Note that we are allowed to use only vertical bridges.

Let $R= [-\frac{\pi}{2}, +\frac{\pi}{2}]\times [0,K_0]$ be a tall rectangle in $\Si$ where $K_0$ be as in Theorem \ref{bridge}. Let $\Sigma_1$ be the unique area minimizing surface with $\PI \Sigma_1=\partial R$. Clearly, $\Sigma_1\simeq S_1$.

We  define $\Sigma_n$ inductively as follows. We will add only vertical bridges to $\Gamma_{n}$ so that the resulting curve $\Gamma_{n+1}$ bounds a unique area minimizing surface $\Sigma_{n+1}$ by Theorem \ref{bridge}.

By Remark \ref{chi}, adding one bridge $\beta_{n+1}$ to $\Sigma_n$ where the endpoints of $\beta_{n+1}$ are in the same component of $\Gamma_n=\PI\Sigma_n$ would suffice to increase the number of ends of $\Sigma_n$ by one. This operation corresponds to {\em adding a pair of pants to the surface}. Similarly by  Remark \ref{chi}, adding two bridges successively so that the endpoints of the first bridge are in the same component, and the endpoints of the second bridge are in different components (components containing the opposite sides of the first bridge),  increases the genus, and keep the number of the ends same. This operation  corresponds to {\em adding a cylinder with handle to the surface}.

Now, we continue inductively to construct the sequence $\{\Sigma_n\}$ dictated by the simple exhaustion (See Figure \ref{simpleexhaustion}). There are two cases: $S_{n+1}-S_n$ contains a pair of pants, or a cylinder with handle.

\vspace{.2cm}

\noindent {\em Pair of pants case.} Assume that $S_{n+1}-S_n$ contains a pair of pants. Let the pair of pants attached to the component $\gamma$ in $\partial S_n$. Let $\gamma'$ be the corresponding component of $\Gamma_n =\PI\Sigma_n$. By construction, $\gamma'$ bounds a disk $D$ in $\Si$ with $D\cap\Gamma_n = \gamma'$. Let $\beta_n=\{c_n\}\times [0,K_0]$ be a vertical segment where $\beta_n\subset D$. Since $\Sigma_n$ bounds a unique area minimizing surface by construction, and $\beta_n\perp \Gamma_n$, we can apply Theorem \ref{bridge}, and get an area minimizing surface $\Sigma_{n+1}$ where $\Sigma_{n+1}$ is homeomorphic to $S_{n+1}$. $\Box$

\vspace{.2cm}

\noindent {\em Cylinder with handle case.} Assume that $S_{n+1}-S_n$ contains a cylinder with handle. Again, let the pair of pants attached to the component $\gamma$ in $\partial S_n$. Let $\gamma'$ be the corresponding component of $\Gamma_n =\PI\Sigma_n$. By construction, $\gamma'$ bounds a disk $D$ in $\Si$ with $D\cap\Gamma_n = \gamma'$. Let $\beta_n$ be a vertical segment $\{c_n\}\times [0,K_0]$ such that $(c_n-\e_n,c_n+\e_n)\times \BR \cap \Gamma_n \subset D$ for some $\e_n>0$. Again, we apply Theorem \ref{bridge} for $\beta_n$ and $\Sigma_n$, and get an area minimizing surface $\Sigma_{n+1}'$. Say $\Gamma_{n+1}'=\PI \Sigma_{n+1}'$. We can choose the thickness of the bridge along $\beta_n$ as small as we want. So, we can assume that the thickness of the bridge along $\beta_n$ is smaller than $\dfrac{\rho.\e_n}{4}$ where $\rho>0$ is the constant in Lemma \ref{disjrec}.

\begin{figure}[b]

\relabelbox  {\epsfxsize=3in

\centerline{\epsfbox{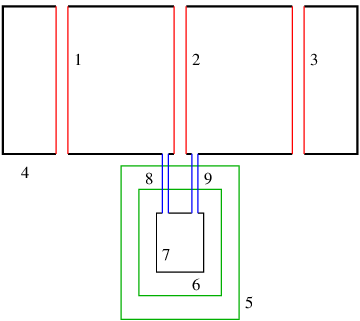}}}

\relabel{1}{\color{red} $\beta_1$}
\relabel{2}{\color{red} $\beta_2$}
\relabel{3}{\color{red} $\beta_3$}
\relabel{4}{$\Gamma_4$}
\relabel{5}{\tiny \color{green} $\partial W_2^+$}
\relabel{6}{\tiny \color{green} $\partial W_2^-$}
\relabel{7}{\tiny $\partial Q_2$}
\relabel{8}{\tiny \color{blue} $\tau_2^-$}
\relabel{9}{\tiny \color{blue} $\tau_2^+$}

\endrelabelbox

\caption{\label{hanger} \footnotesize In the figure above, $S_2-S_1$ is a pair of pants, and $S_3-S_2$ is a cylinder with handle. When constructing $\Sigma_3$, $\beta_2$ is attached to the corresponding component in $\Gamma_2$, then a hanger, the pair of vertical bridges $\tau^\pm_2$ and a thin rectangle $Q_n$, is added to obtain the cylinder with handle. $\partial W_2^\pm$ is needed to show that $\Sigma_2' \cup T_2$ is uniquely area minimizing surface to apply Theorem \ref{bridge}. }

\end{figure}

Now, consider the rectangle $Q_n = [c_n-\frac{\rho.\e_n}{2}, c_n+\frac{\rho.\e_n}{2}]\times [-6\pi-K_0,-4\pi-K_0]$ (See Figure \ref{hanger}). Let $T_n$ be the unique area minimizing surface in $\BHH$ with $\PI T_n=\partial Q_n$ by Lemma \ref{rectangle}. Let $\wh{\Gamma}_{n+1}=\Gamma_{n+1}'\cup \partial Q_n$. We claim that $\wh{\Gamma}_{n+1}$ bounds a unique area minimizing surface $\wh{\Sigma}_{n+1}$ in $\BHH$ and $\wh{\Sigma}_{n+1}=\Sigma_{n+1}'\cup T_n$. Notice that $\Sigma_{n+1}'$ and $T_n$ are uniquely minimizing surfaces. Hence, if we show that $\wh{\Gamma}_{n+1}$ cannot bound any connected area minimizing surface, then we are done.

Assume that $\wh{\Gamma}_{n+1}$ bounds a connected area minimizing surface $\wh{\Sigma}_{n+1}$. Consider the the pair of rectangles $W^+_n=[c_n-\e_n,c_n+\e_n]\times [-9\pi-K_0,-\pi-K_0] $ and $W^-_n=[c_n-\rho.\e_n,c_n+\rho.\e_n]\times [-7\pi-K_0,-3\pi-K_0]$. Let $\Upsilon_n=\partial W_n^+\cup\partial W_n^-$. Then, by Lemma \ref{disjrec}, the uniquely minimizing surface $F_n$ with $\PI F_n= \Upsilon_n$ must be $P^+_n\cup P^-_n$ where $P^\pm_n$ is the unique area minimizing surface with $\PI P^\pm_n = \partial W^\pm_n$. As $\wh{\Gamma}_{n+1}\cap \Upsilon_n =\emptyset$, the area minimizing surfaces $\wh{\Gamma}_{n+1}$ and $F_n$ must be disjoint by Lemma \ref{disjoint} (See Figure \ref{hanger}). On the other hand, the area minimizing surface $F_n=P^+_n\cup P^-_n$ separates the components, $\Gamma_{n+1}'$  and $\partial Q_n$, of $\wh{\Gamma}_{n+1}$. Since $\wh{\Gamma}_{n+1}\cap F_n=\emptyset$, this implies $\wh{\Sigma}_{n+1}$ disconnected. This proves that $\wh{\Sigma}_{n+1}=\Sigma_{n+1}'\cup T_n$ is the unique area minimizing surface with $\PI \wh{\Sigma}_{n+1}= \wh{\Gamma}_{n+1}$.

Now, let $\tau^+_n=\{c_n+\frac{\rho.\e_n}{4}\}\times [-4\pi-K_0,0]$ be the vertical arc segment in $\Si$. When we apply Theorem \ref{bridge} to the uniquely minimizing surface $\wh{\Sigma}_{n+1}$ and the arc $\tau^+_n$, we obtain a  new uniquely minimizing surface $\wh{\Sigma}'_{n+1}$. Similarly, let $\tau^-_n=\{c_n-\frac{\rho.\e_n}{4}\}\times [-4\pi-K_0,0]$. Again, we apply Theorem \ref{bridge} for $\wh{\Sigma}'_{n+1}$ and $\tau^-_n$, we obtain another uniquely minimizing surface $\Sigma_{n+1}$. Furthermore, we  assume that the both bridges along $\tau^+_n$ and $\tau^-_n$ have thickness less than $\frac{\rho.\e_n}{4}$. The pair of vertical bridges along $\tau^\pm_n$ with the thin rectangle $Q_n$ looks like a hanging picture frame (See Figure \ref{hanger}).

By construction, $\Sigma_{n+1}$ is homeomorphic to $S_{n+1}$. In particular, we achieved to add a cylinder with handle to $\Sigma_n$ along the corresponding component $\gamma'$ in $\Gamma_n$. This finishes the description of the inductive step, when $S_{n+1}-S_n$ contains a cylinder with handle. $\Box$\\

\noindent {\em The Limit and the Properly Embeddedness:} Notice that in the bridge principle at infinity (Theorem \ref{bridge}), as the thickness of the bridge $\alpha$ goes to $0$, the height of the strip $S_\alpha$ goes to $0$, too. In particular, let $\Gamma,\Sigma, \alpha, \Gamma_t,\Sigma_t$ be as in the statement of Theorem \ref{bridge}. Let $S^t_\alpha =\Sigma_t\cap N_\e(\alpha)$ where $N_\e(\alpha)$ is the sufficiently small neighborhood of $\alpha$ in the compactification $\overline{\BHH}$. Then, as $t\searrow0$, then $d(L_z, S^t_\alpha)\to \infty$ where $L_z$ is the vertical line through origin in $\BHH$, i.e. $L_z=\{0\}\times \BR$. This is because as $t\searrow 0$, $\Sigma_t\to \Sigma$.

Let $\wh{B}_r=B_r(0)\times [-2K_0,2K_0]$ be compact region in $\BHH$ where $B_r(0)$ is the $r$ ball around origin in $\BH^2$. As $t_n\searrow 0$, then the thickness of the bridge in $\Sigma_n$ near $\beta_n$ (or $\tau_n^\pm$) goes to $0$. Hence, by choosing $t_n<\frac{1}{10n^2}$ sufficiently small, we can make sure that $d(L_z, S^{t_n}_{\beta_n})>r_n$ and $d(L_z, S^{t_n}_{\tau_n^\pm})>r_n$ for a sequence $r_n\nearrow \infty$. This implies that for $m\geq n$, $\wh{B}_{r_n}\cap \Sigma_m \simeq S_n$, as the thickness (and hence height) of the bridges $\beta_n$ and $\zeta_n$ goes to $0$.

Now, $\Sigma_n$ is a sequence of absolutely area minimizing surfaces in $\BHH$. Let $\Sigma_n'=\wh{B}_{r_n}\cap\Sigma_n$. By Lemma \ref{convergence}, by using a diagonal sequence argument, we get a limiting surface $\Sigma$ in $\BHH$ where the convergence is smooth on compact sets. $\Sigma$ is an area minimizing surface in $\BHH$ as it is the limit of area minimizing surfaces. Notice that for $m\geq n$, $\wh{B}_{r_n}\cap \Sigma_m \simeq S_n$ and the convergence is smooth on compact sets. This implies $\Sigma\cap \wh{B}_{r_n}\simeq S_n$ for any $n$, and hence $\Sigma\simeq S$.

We also note that the bridges do not collapse in the limit, as for every bridge along $\beta_n$ and $\tau_n^\pm$, we can place a thin, tall rectangle $R_n$ "under" the bridge {\em disjoint} from the minimizing sequence. In other words, the area minimizing plane $P_n$ with $\PI P_n=R_n$ (Lemma \ref{rectangle}) will be a barrier for bridges to collapse, as  for any $m>n$, $P_n\cap\Sigma_m=\emptyset$ since $\Gamma_m\cap R_n=\emptyset$ by Lemma \ref{disjoint}.

Finally, $\Sigma$ is properly embedded in $\BHH$ as for any compact set $K\subset \BHH$, there exists $r_n>0$ with $K\subset \wh{B}_{r_n}$, and $\wh{B}_{r_n}\cap \Sigma\simeq S_n$ which is compact. The proof of the theorem follows.
\end{pf}

\section{Final Remarks} \label{secremarks}

\subsection{H-surfaces} \

\vspace{.2cm}

The constant mean curvature surfaces could be considered as a natural candidate to generalize our results. Hence, consider the following question:

\vspace{.2cm}

{\bf Question:} {\em What kind of surfaces can be embedded in $\BHH$ as a complete $H$-surface for $0<H<\frac{1}{2}$?}

\vspace{.2cm}

In other words, is it possible to embed any open orientable surface $S$ in $\BHH$ as a complete $H$-surface for $0<H<\frac{1}{2}$. A positive answer to these question would be a generalization of Theorem \ref{main} to $H$-surfaces.

Unfortunately, it is hardly possible to generalize our methods to this problem. By \cite{NSST}, for $H>0$, if $\Sigma$ is an $H$-surface with $\PI\Sigma \neq \emptyset$ and $\Sigma\cup\PI\Sigma$ is a $C^1$ surface up to the boundary, then $\PI\Sigma$ must be a collection of a vertical line segments in $\Si$. In particular, this implies the asymptotic Plateau problem practically has no solution for $H$-surfaces in $\BHH$ since if $\Gamma$ is a $C^1$ simple closed curve in $\Si$, there is no $H$-surface $\Sigma$ in $\BHH$ where $\Sigma\cup \Gamma$ is a $C^1$ surface up to the boundary. Hence, because of this result, our methods for Theorem \ref{main} cannot be generalized to this case. However, it might be possible to construct a complete $H$-surface $\Sigma$ of any finite topology with only vertical ends, i.e. $\PI \Sigma$ consists of only vertical lines in $\Si$.


\subsection{Finite Total Curvature.} \

\vspace{.2cm}

Our construction of area minimizing surfaces In $\BHH$ produces surfaces of infinite total curvature. In \cite{MMR}, Martin, Mazzeo and Rodriguez recently showed that for any $g\geq 0$, there exists a complete, finite total curvature, embedded minimal surface $\Sigma_{g,k_g}$ in $\BHH$ with genus $g$ and $k_g$ ends for sufficiently large $k_g$. Even though this result is a great progress to construct examples of minimal surfaces of finite total curvature, the question of existence (or nonexistence) of minimal surfaces of finite total curvature with any finite topology is still a very interesting open problem.

It is well known that a complete, properly embedded, minimal surface in $\BHH$ with finite total curvature has also finite topology \cite{HR}. On the other hand, there are surfaces with finite topology which cannot be embedded in $\BHH$ as a complete minimal surface with finite total curvature. For example, by \cite{HNST}, a twice punctured torus cannot be embedded as a complete minimal surface with finite total curvature into $\BHH$. Hence, the following question becomes very interesting:\\

\noindent {\bf Question:} {\em For which $g\geq0$, and $k\geq 0$, there exists a complete embedded minimal surface $S^g_k$ in $\BHH$ with finite total curvature where $S^g_k$ is an orientable surface of genus $g$ with $k$ ends?}

\section{Appendix} \label{secappendix}

In this section, we prove some technical steps used in our construction.

\subsection{Generic Uniqueness of Area Minimizing Surfaces} \label{GenUniqSec} \

\vspace{.2cm}

In this part, we prove a generic uniqueness result for tall curves in $\Si$. Note that the results in this part are mostly for area minimizing surfaces, and do not apply to minimal surfaces in general.

We start with a lemma which roughly says that disjoint curves in $\Si$ bounds disjoint area minimizing surfaces in $\BHH$.

\begin{lem} [Disjointness] \label{disjoint} Let $\Omega_1$ and $\Omega_2$ be two closed regions in $\SI$ where $\partial \Omega_i=\Gamma_i$ is a finite collection of disjoint simple closed curves. Further assume that $\Omega_1\cap\Omega_2=\emptyset$ or $\Omega_1\subset int(\Omega_2)$. If $\Sigma_1$ and $\Sigma_2$ are area minimizing surfaces in $\BHH$ with $\PI\Sigma_i = \Gamma_i$, then $\Sigma_1\cap\Sigma_2 =\emptyset$.
\end{lem}

\begin{pf} Assume that $\Sigma_1\cap\Sigma_2\neq \emptyset$. As both surfaces are minimal, by maximum principle, the intersection cannot contain isolated points.  As $\Gamma_1\cap \Gamma_2=\emptyset$, then $\Sigma_1\cap\Sigma_2=\alpha$ which is collection of closed curves.
	
	Since $\BHH$ is topologically a ball, any surface would be separating. Let $\Delta_i$ be the components of $\BHH-\Sigma_i$ with $\PI \overline{\Delta}_i = \Omega_i$. In other words, as $\Sigma_i\cup\Omega_i$ is a closed surface in the contractible space $\overline{\BHH}$, it bounds a region $\Delta_i$ in $\overline{\BHH}$.
	
	If $\Omega_1\subset int(\Omega_2)$, let $S_1=\Sigma_1-\Delta_2$ and let $S_2=\Sigma_2\cap \overline{\Delta}_1$. Then, as $\Omega_1\subset int(\Omega_2)$, with this operation, we cut the surfaces $S_i$ from the non-compact parts in $\Sigma_i$. Therefore, $\PI S_1=\PI S_2=\emptyset$ and both $S_1$ and $S_2$ are compact surfaces with $\partial S_1=\partial S_2 = \alpha$.
	
	If $\Omega_1\cap\Omega_2=\emptyset$, let $S_1=\Sigma_1\cap\overline{\Delta}_2$ and let $S_2=\Sigma_2\cap \overline{\Delta}_1$. Again, as $\Omega_1\cap\Omega_2=\emptyset$, $\PI S_1=\PI S_2=\emptyset$ and both $S_1$ and $S_2$ are compact surfaces with $\partial S_1=\partial S_2 = \alpha$.
	
	As $\Sigma_1$ and $\Sigma_2$ are area minimizing surfaces, so are $S_1\subset \Sigma_1$ and $S_2\subset \Sigma_2$. Hence, as $\partial S_1=\partial S_2$, $|S_1|=|S_2|$ where $|.|$ represents the area. Let $T_1$ be a compact subsurface in $\Sigma_1$ containing $S_1$, i.e. $S_1\subset T_1\subset \Sigma_1$. Consider $T_1'=(T_1-S_1)\cup S_2$. Since $T_1$ is area minimizing and $|T_1'|=|T_1|$, so is $T_1'$. However, $T_1'$ is not smooth along $\alpha$ which contradicts to the regularity of area minimizing surfaces (Lemma \ref{AMSexist}). The proof follows.
\end{pf}

\begin{rmk} \label{disjointrem} Note that in the lemma above $\Omega_1\cap\Omega_2=\emptyset$ or $\Omega_1\subset int(\Omega_2)$ are indeed equivalent conditions. This is because we can always replace $\Omega_2$ with $\overline{\Omega_2^c}$. Notice also that the proof above is simply a swaping argument ($S_1$ and $S_2$) for area minimizing surfaces, and the proof actually works for more general case. In particular, we do not need $\Gamma_i$ to be a collection of {\em simple} closed curves, but only to be $\Gamma_i=\partial \Omega_i$ where $\Omega_1\subset int(\Omega_2)$ for swaping argument. So, with the same proof, the lemma above can also be stated as follows: {\em Let $\Omega_1$ and $\Omega_2$ be two open regions in $\SI$ with $\overline{\Omega}_1\subset \Omega_2$. If $\Sigma_1$ and $\Sigma_2$ are area minimizing surfaces in $\BHH$ with $\PI\Sigma_i = \partial\overline{\Omega}_i$, then $\Sigma_1\cap\Sigma_2 =\emptyset$.}
\end{rmk}

Now, we show that if a tall curve $\Gamma\subset\Si$ does not bound a unique area minimizing surface in $\BHH$, it bounds two canonical area minimizing surfaces $\Sigma^\pm$ where any other area minimizing surface $\Sigma'$ with $\PI\Sigma'=\Gamma$ must be "between" $\Sigma^+$ and $\Sigma^-$.

\begin{lem} [Canonical Surfaces] \label{canonical}
	Let $\Gamma$ be a tall curve in $\Si$. Then either there exists a unique area minimizing
	surface $\Sigma$ in $\BHH$ with $\PI \Sigma = \Gamma$, or there are two canonical disjoint extremal
	area minimizing surfaces $\Sigma^+$ and $\Sigma^-$ in $\BHH$ with $\PI \Sigma^\pm = \Gamma$.
\end{lem}

\begin{pf} We  mainly adapt the techniques of \cite[Lemma 4.3]{Co2} (Similar result for $\BH^3$) to $\BHH$ context. Let $\Gamma$ be a tall curve in $\Si$, and let $\Gamma^c=\Omega^+\cup \Omega^-$ where $\Omega^\pm$ are two tall regions in $\Si$ with $\partial\overline{\Omega^+}=\partial\overline{\Omega^-}=\Gamma$. Let $N_\e(\Gamma)$ be a small open neighborhood of $\Gamma$ in $\Si$. Let $N^+=N_\e(\Gamma)\cap \Omega^+$ and let $N^-=N_\e(\Gamma)\cap \Omega^-$. Let the family of curves $\{\Gamma^\pm_t\mid t\in [0,\e)\}$ foliate the region $N^\pm$ with $\Gamma_0=\Gamma$. Let $\Gamma^\pm_n=\Gamma^\pm_{t_n}$ for $t_n\searrow 0$. By choosing $\e>0$ sufficiently small, we can assume $\Gamma^\pm_n$ is tall for any $n>0$. Let $\Sigma^\pm_n$ be an area minimizing surface in $\BHH$ with $\PI \Sigma^\pm_n=\Gamma^\pm_n$ by Lemma \ref{APP}.
	
	By replacing the sequence $\Sigma_n$ with $\wh{B}_n\cap\Sigma^\pm_n$ in the proof of Lemma \ref{APP}, we can show that $\Sigma^+_n$ converges (up to a subsequence) to an area minimizing surface $\Sigma^+$ with $\PI \Sigma^+ = \Gamma$. Similarly, $\Sigma^-_n$ converges to an area minimizing surface $\Sigma^-$ with $\PI \Sigma^- = \Gamma$.
	
	Assume that $\Sigma^+ \neq\Sigma^-$, and they are not disjoint. By maximum principle, they cannot have isolated points in the intersection. Therefore, nontrivial intersection implies some part of $\Sigma^-$ lies {\em above} $\Sigma^+$, i.e. some part of $\Sigma^-$ separated by $\Sigma^=$. Then, since	$\Sigma^+=\lim \Sigma_n^+$, $\Sigma^-$ must also intersect some $\Sigma_n^+$ for sufficiently large $n$. However by Lemma \ref{disjoint} (Swaping argument), $\Sigma_n^+$ is disjoint from $\Sigma^-$ as $\PI\Sigma_n^+ = \Gamma_n^+$ is disjoint from $\PI \Sigma^- = \Gamma$. This is a contradiction. This shows $\Sigma^+$ and $\Sigma^-$ are disjoint. By using similar techniques to \cite[Lemma 4.3]{Co2}, it can be showed that $\Sigma^\pm$ are canonical, i.e. independent of the sequences $\{\Sigma^\pm_n\}$.
	
	Similar arguments show that $\Sigma^\pm$ are disjoint from any area minimizing hypersurface $\Sigma '$	with $\PI\Sigma'=\Gamma$. As the sequences of $\Sigma_n^+$ and $\Sigma_n^-$ forms a barrier for other area	minimizing hypersurfaces asymptotic to $\Gamma$, any such area minimizing hypersurface must lie in the	region bounded by $\Sigma^+$ and $\Sigma^-$ in $\BHH$. This shows that if $\Sigma^+ = \Sigma^-$, then there exists a	unique area minimizing hypersurface asymptotic to $\Gamma$. The proof follows.
\end{pf}

\begin{rmk} Notice that if a finite collection of simple closed curves $\Gamma$ is not assumed to be tall in the lemma above, the same proof is still valid. Hence, for any such $\Gamma$, either there is either no solution ($\nexists \Sigma$), or a unique solution ($\exists ! \Sigma$), or two canonical solutions ($\exists \Sigma^\pm$) for asymptotic Plateau problem for $\Gamma$ ($\PI\Sigma = \Gamma$).
\end{rmk}

Now, by using the lemma above, we show a generic uniqueness result for tall curves.

\begin{thm} [Generic Uniqueness] \label{uniq} A generic tall curve in $\Si$ bounds a unique area minimizing surface in $\BHH$.
\end{thm}

\begin{pf} Let $\Gamma_0$ be a tall curve in $\Si$. Let $N(\Gamma_0)$ be a small open neighborhood of $\Gamma_0$ in $\Si$ which is a finite collection of annuli. Let $\{\Gamma_t \ | \ t\in(-\e,\e)\}$ be a foliation of $N(\Gamma_0)$. In particular, for any $-\e<t_1<t_2<\e$, $\Gamma_{t_1}\cap \Gamma_{t_2}=\emptyset$. We can assume $N(\Gamma_0)$ sufficiently thin that $\Gamma_t$ is a tall curve for any $t\in (-\e,\e)$. Let $\Sigma_t$ be an area minimizing surface in $\BHH$ with $\PI \Sigma_t=\Gamma_t$.
	
	As in the proof of the lemma above, let $\Gamma^c_t=\Omega_t^+\cup \Omega_t^-$ with $\partial\overline{\Omega_t^+}=\partial\overline{\Omega_t^-}=\Gamma_t$. Then, $\Omega^+_t\subset\Omega^+_s$ for $t<s$. Hence by Lemma \ref{disjoint}, $\Sigma_{t}\cap\Sigma_{s}=\emptyset$ for $t<s$. Furthermore, by Lemma \ref{canonical}, if $\Gamma_s$ does not bound a unique area minimizing surface $\Sigma_s$, then we can define two disjoint canonical minimizing  $\Sigma^+_s$ and $\Sigma^-_s$ with $\PI \Sigma^\pm_s = \Gamma_s$. Hence, $\Sigma_s^+\cup \Sigma_s^-$ separates a region $V_s$ in $\BHH$. If $\Gamma_s$ bounds a unique area minimizing surface $\Sigma_s$, then let $V_s=\Sigma_s$ (say $V_s$ a degenerate neighborhood). Notice that by lemma \ref{disjoint}, $\Sigma_t\cap \Sigma_s=\emptyset$ for $t\neq s$, and hence $V_t\cap V_s=\emptyset$ for $t\neq s$.
	
	
	Now, consider a short arc segment $\eta$ in $\BHH$ with one endpoint is in $\Sigma_{t_1}$ and the other end point is in $\Sigma_{t_2}$ where $-\e<t_1<0<t_2<-\e$. Hence, $\eta$ intersects all area minimizing surfaces $\Sigma_t$ with $\PI \Sigma_t = \Gamma_t$ where $t_1\leq t \leq t_2$. Now for $t_1<s<t_2$, define the {\em thickness} $\lambda_s$ of $V_s$ as $\lambda_s=|\eta\cap V_s|$, i.e. $\lambda_s$ is the length of the piece of $\eta$ in $V_s$. Hence, if $\Gamma_s$ bounds more than one area minimizing surface, then the thickness $\lambda_s>0$. In other words, if $\lambda_s=0$, then $\Gamma_s$ bounds a unique area minimizing surface in $\BHH$.
	
	As $V_t\cap V_s=\emptyset$ for $t\neq s$, we have $\sum_{t_1}^{t_2}\lambda_s<|\eta|$. Hence, as $|\eta|$ is finite, $\lambda_s>0$ for only countably many $s\in[t_1,t_2]$. This implies for all but countably many $s\in[t_1,t_2]$, $\lambda_s=0$, and hence $\Gamma_s$ bounds a unique area minimizing surface. Similarly, this implies for all but countably many $s\in(-\e,\e)$, $\Gamma_s$ bounds a unique area minimizing surface. Then, by using the techniques in \cite[Lemma 3.2]{Co2}, the generic uniqueness in Baire Sense follows.
\end{pf}

\subsection{Nonexistence Results for Vertical Bridge Principle:} \

\vspace{.2cm}

The following lemma rules out some special cases for asymptotic Plateau problem, and used in the proof of the bridge principle. 

\begin{lem} \label{circles} Let $\gamma_c=\si\times\{c\}$ represent the round circle in $\Si$ with $\{z=c\}$.  Let $\Gamma=\bigcup_{i=1}^N\gamma_{c_i}\bigcup_{j=1}^M \alpha_j$ where $\alpha_j=\{\theta_j\}\times[c_{j_1},c_{j_2}]$ for some $\theta_j\in\si$, and $c_i<c_{i+1}$. Then, there exists a $K_0>\pi$ such that the following holds: If $c_{i+1}-c_i>K_0$ for any $i$,	and $\Sigma$ is an area minimizing surface in $\BHH$ with $\PI \Sigma\subset\Gamma$, then $\Sigma$ is a collection of horizontal planes, i.e. $M=0$ and $\Sigma=\BH^2\times\{c_{i_1},c_{i_2},...c_{i_k}\}$.
\end{lem}

Notice that the statement implies for such $K_0>0$, the asymptotic boundary of such an area minimizing surface cannot contain any vertical line segment.

\begin{pf} Without loss of generality, we  assume $N=2$ as the other cases are similar. We  divide the proof into two cases: $M=1$ and $M>1$.

\vspace{.2cm}

\noindent {\bf Case 1:} Assume $M=1$, i.e. $\Gamma=\gamma_{c_1}\cup\gamma_{c_2}\cup\alpha_1$ where $\alpha_1=\{\theta_1\}\times[c_1,c_2]$ for some $\theta_1\in\si$. Let $\Sigma$ be the area minimizing surface in $\BHH$ with $\PI \Sigma\subset\Gamma$. Recall that $c_2-c_1>K_0>\pi$. Let $R_i$ be a sequence of tall rectangles exhausting the region bounded by $\Gamma$, i.e.  $R_i=\partial ( [\theta_1+\e_i,\theta_1-\e_i+2\pi]\times [c_1+\rho_i,c_2-\rho_i])$ in $\Si$ where $\e_i\searrow 0$ and $\rho_i\searrow 0$. Clearly, $R_i$ is disjoint from $\Gamma$ for any $i$, and $R_i\to \Gamma$ as $i\to \infty$. 
	
Let $P_i$ be the unique area minimizing surface in $\BHH$ with $\PI P_i = R_i$ (Lemma \ref{rectangle}). By Lemma \ref{disjoint} and Remark \ref{disjointrem}, $\Sigma\cap P_i=\emptyset$ for any $i$.	On the other hand, the explicit description of $P_i$ in \cite{ST} shows that $P_i$ is foliated by horizontal equidistant curves $\beta_i^t =P_i\cap \BH^2\times\{t\}$ to the geodesic $\tau_i$ with  $\PI\tau_i=\{\theta_1+\e_i,\theta_1-\e_i+2\pi\}$.  In particular, for $d_i(t)=d(\beta_i^t,\tau_i)$,  $d_i(t)\to \infty$ as $t\to c_1$ or $t\to c_2$, while $d_i(c*)<C_0$ where $c*=\frac{c_1+c_2}{2}$ (See the discussion before Lemma \ref{rectangle}).	Hence, as $i\to \infty$, $\tau_i$ and hence $\beta_i^{c*}$ escapes to infinity. This shows $P_i$ converges to two horizontal geodesic planes $\BH^2\times\{c_1,c_2\}$. However, this implies $\Sigma\cap P_i\neq \emptyset$ for sufficiently large $i$ unless $\PI \Sigma\subset \gamma_{c_1}\cup\gamma_{c_2}$. Hence, $M=1$ case follows. \hfill $\Box$

\vspace{.2cm}

\noindent {\bf Case 2:} Now, assume $M>1$. By using a simple trick, we  reduce this case to $M=2$. Let $\theta_0\in \si-\{\theta_1,\theta_2,...\theta_M\}$. Let $\tau$ be the geodesic in $\BH^2$ with $\PI \tau = \{\theta_0,\theta_1\}$. Let $\varphi$ be the hyperbolic isometry fixing $\tau$ pushing from $\theta_1$ to $\theta_0$ with translation length $l>0$. Let $\wh{\varphi}$ be the isometry of $\BHH$ with $\wh{\varphi}(x,t)=(\varphi(x),t)$. Then, define the sequence of area minimizing surfaces $\Sigma_n=\wh{\varphi}^n(\Sigma)$. Then, by Lemma \ref{convergence}, there exists a subsequence of $\{\Sigma_n\}$ converging to an area minimizing surface $\wh{\Sigma}$ in $\BHH$. Let $\wh{\Gamma}=\PI\wh{\Sigma}$. By construction, $\wh{\Sigma}$ is invariant under $\wh{\varphi}$, then so is $\wh{\Gamma}$. As $\{\theta_0,\theta_1\}$ are the fixed points of $\varphi$, this implies $\wh{\Gamma}\subset \gamma_{c_1}\cup\gamma_{c_2}\cup\alpha_0\cup \alpha_1$ where $\alpha_i=\{\theta_i\}\times[c_1,c_2]$. 

We claim that $\wh{\Gamma} = \gamma_{c_1}\cup\gamma_{c_2}\cup\alpha_0\cup \alpha_1$. Clearly, $\wh{\Gamma}\supset \gamma_{c_1}\cup\gamma_{c_2}$ by construction. Now, consider a component $S$ of $\Sigma$ with $\PI S \supset \alpha_1$ (possibly $\Sigma=S$). Since, we assumed $M>1$, $\PI S$ must contain another $\alpha_{j_0}$ for some $j_0>1$. By  Lemma \ref{surface}, $\overline{S}=S\cup\PI S$ is a surface with boundary in $\overline{\BHH}$.
Consider the collection of curves $\lambda_c=S \cap \BH^2\times \{c\}$ for $c\in (c_1,c_2)$. As $S$ is connected, there exist a $c\in (c_1,c_2)$ such that $\lambda_c$ contains a infinite line $l_c$ in $\BH^2\times \{c\}$ with $\PI l_c=\{\theta_1,\theta_j\}$. Let $l^n_c=\wh{\varphi}^n(l_c)\subset \Sigma_n\cap \BH^2\times \{c\}$. Then, by construction $l^n_c$ converges to $\wh{l}_c\subset \wh{\Sigma}\cap\BH^2\times \{c\}$ where $\PI \wh{l}_c = \{\theta_1,\theta_0\}$. This proves that $\alpha_0\cup \alpha_1 \supset \wh{\Gamma}$. Hence, we reduce the $M>1$ case to $M=2$ case.

Now, we  finish this case. Recall that by construction $\wh{\Sigma}$ is invariant by $\wh{\Sigma}$, i.e. $\wh{\varphi}(\wh{\Sigma})=\wh{\Sigma}$. Because of this invariance, we first claim that $\wh{\Sigma}=\p_0\cup\p_1$ where $\p_i$ is the unique area minimizing plane with asymptotic boundary a rectangle $R_i$, i.e. $\PI\p_0 =R_0=\partial ([\theta_0,\theta_1]\times [c_1,c_2])$ and $\PI\p_1 =R_1=\partial ([\theta_1,\theta_0+2\pi]\times [c_1,c_2])$. In order to see this, 
let $\theta_2=\frac{\theta_0+\theta_1}{2}$ and $\theta_3=\theta_2+\pi$ in $\si$. Let $\eta$ be the geodesic in $\BH^2$ with $\PI\eta=\{\theta_2,\theta_3\}$. Let $\mathcal{W}=\eta\times \BR$ be the vertical plane in $\BHH$. Consider 
$Z=\mathcal{W}\cap\wh{\Sigma}$. By construction, $Z$ is a collection of curves with $\PI Z$ is the four points, $(\theta_2,c_1), (\theta_2,c_2), (\theta_3,c_1), (\theta_3,c_2)$. Invariance of $\wh{\Sigma}$ by $\wh{\varphi}$ implies that $Z$ is the generating curves for $\wh{\Sigma}$. Assuming $\wh{\Sigma}\neq \BH^2\times \{c_1,c_2\}$, by \cite{ST}, we conclude that $Z=\mu_0\cup \mu_1$ where $\mu_0$ is the generating curve for $\p_0$, and $\mu_1$ is the generating curve for $\p_1$ such that $\PI \mu_0= \{(\theta_2,c_1), (\theta_2,c_2)\}$ and $\PI \mu_1= \{(\theta_3,c_1), (\theta_3,c_2)\}$. Now, even though the union $\p_0\cup\p_1$ is a minimal surface in $\BHH$, we  show that it is not an area minimizing surface, and finish the proof of Case 2. 

\vspace{.2cm}

\noindent {\bf Claim:} $\p_0\cup\p_1$ is not an area minimizing surface.

\vspace{.2cm}

We  show that a sufficiently long annulus $\A$ between $\p_0$ and $\p_1$ has less area than the sum of the areas of the corresponding disks $D_0$ in $\p_0$ and $D_1$ in $\p_1$, i.e. $\partial \A= \partial D_0\cup\partial D_1$ (See Figure \ref{banana}).

Without loss of generality, let $c_1=-K$ and $c_2=K$, and $\theta_0=0$ and $\theta_1=\pi$ in $\si$. By \cite[Proposition 2.1 (1)]{ST} and Lemma \ref{rectangle}, we have a very good understanding of the area minimizing planes $\p_0$ and $\p_1$. By the symmetry, we  work with only $\p_0$. Let $\upsilon$ be the geodesic in $\BH^2$ with $\PI \upsilon=\{\frac{\pi}{2},\frac{3\pi}{2}\}$.
Recall that $\p_0$ has the generating curve $c_0$ in the vertical plane $\upsilon\times\BR$ where $\PI c_0 = \{(\frac{\pi}{2}, -K), (\frac{\pi}{2}, K)$. The parametrization of the generating curve $c_0$ has explicitly been given in the proof of \cite[Proposition 2.1]{ST} as $\lambda(\rho)$ for $d>1$ (case (1)). Now, recall that $\PI \p_0 = R_0= \partial ([0,\pi]\times[-K,K])$.  Let $t\in [-K,K]$ represent the height in $\BHH$. Parametrize $c_0$ as $\lambda(t)=(\rho(t), t)$ in $\upsilon\times\BR$, where $\rho(t)$ is the distance of $(0,t)$ to $c_0\cap \BH^2\times \{t\}$. 

Recall that $\tau$ is the geodesic in $\BH^2$ with $\PI \tau=\{0,\pi\}$. Parametrize $\tau$ such that $\tau(s)$ is the signed distance from the origin for $s\in (-\infty,+\infty)$. In particular, $\tau(+\infty)=\pi$ and $\tau(+\infty)=\pi$ in $\si$. Let $\varphi_t$ be the hyperbolic isometry fixing $\tau$ with translation length $t\in\BR$. Then, by \cite{ST}, $\varphi_t(\p_0)=\p_0$ for any $t\in \BR$.
Let $\p_o\cap \BH^2\times \{t\}= \eta_t$. Then by construction, $\eta_t$ is the equidistant line to $\tau$ with distance $\rho(t)$. Parametrize $\eta_t$ such that the closest point to $\tau(s)$ in $\eta_t$ would be $\eta_t(s)$ for $s\in (-\infty,+\infty)$.

\begin{figure}[t]
	\begin{center}
		$\begin{array}{c@{\hspace{.2in}}c}

		\relabelbox  {\epsfysize=2in \epsfbox{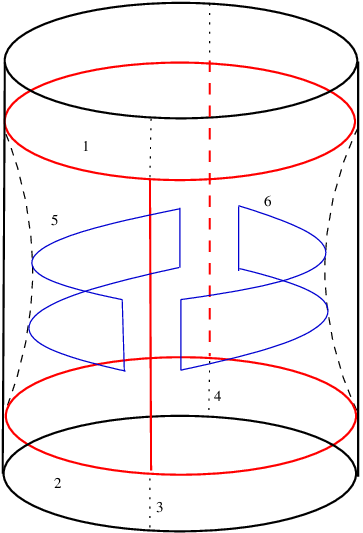}} \relabel{1}{\tiny $\gamma_K$} \relabel{2}{\tiny $\gamma_{-K}$} 
		\relabel{3}{\tiny $\theta=0$} \relabel{4}{\tiny $\theta=\pi$} \relabel{5}{\tiny $\partial D_1$} \relabel{6}{\tiny $\partial D_0$} 		\endrelabelbox & 
		
		\relabelbox  {\epsfysize=2in \epsfbox{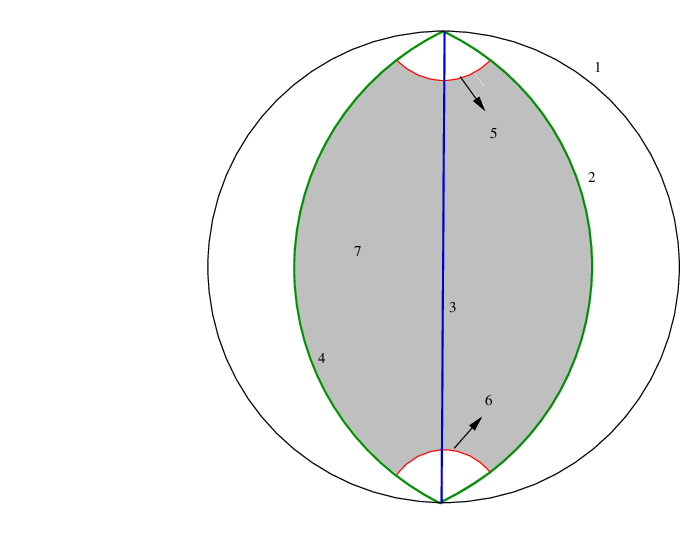}} \relabel{1}{\tiny $\BH^2\times\{k\}$} \relabel{2}{\tiny $\xi_0^+$} 
		\relabel{3}{\tiny $\tau$} \relabel{4}{\tiny $\xi_1^+$} \relabel{5}{\tiny $\sigma^+_+$} \relabel{6}{\tiny $\sigma^+_-$} \relabel{7}{\tiny $\Delta^+$} 
		\endrelabelbox\\
		\end{array}$
		
	\end{center}
	\caption{\label{banana} \footnotesize In the left, red curve represents $\wh{\Gamma}$ in $\Si$. Blue curves represents $\partial D_i$ in $\p_i$. In the right, the domain $\Delta^+$ is depicted in the banana region between equidistant lines to $\tau$.}
\end{figure}

Now, we describe $D_i$ in $\p_i$ by defining its boundary $\partial D_i$.  Like $\p_0$ and $\p_1$, $D_0$ and $D_1$ will be symmetric with respect to $\tc=\tau\times \BR$ so let's only consider $D_0$. $\partial D_0$ is a rectangle in $\p_0$ with the following four edges. Fix $k>\frac{\pi}{2}$ be the half height of rectangle with $k<<K$.
Let the upper edge $\xi_0^+$ be the segment in $\eta_k$ between the points $\eta_k(-l)$ and $\eta_k(l)$ where $l>>0$ will be determined later. Similarly, the lower edge $\xi_0^-$ be the be the segment in $\eta_{-k}$ between the points $\eta_{-k}(-l)$ and $\eta_{-k}(l)$. Let the short edges be the vertical paths $\nu_0^+$ and $\nu_0^-$ in $\p_0$ with endpoints $\{\eta_k(l),\eta_{-k}(l)\}$ and  $\{\eta_k(-l),\eta_{-k}(-l)\}$ respectively. Hence, $D_0$ is the rectangle in $\p_0$ with $\partial D_0= \xi_0^+\cup \nu_0^+\cup\xi_0^-\cup\nu_0^-$. Similarly, define $D_1$ in $\p_1$ as $\partial D_1= \xi_1^+\cup \nu_1^+\cup\xi_1^-\cup\nu_1^-$ as the symmetric rectangle with respect to the vertical plane $\tc$ (See Figure \ref{banana}-left).

Now, we define the {\em competitor} annulus $\A$ with $\partial \A=\partial D_0\cup\partial D_1$. Let $\sigma^+_+$ be the geodesic between  $\eta_k(l)$ and its reflection with respect to $\tc$. Let $\sigma^+_-$ be the reflection of $\sigma^+_+$ with respect to $\upsilon\times\BR$. Let $\sigma^-_+$ be the reflection of $\sigma^+_+$ with respect to horizontal plane $\BH^2\times\{0\}$. Similarly, let $\sigma^-_-$ be the reflection of $\sigma^+_-$ with respect to horizontal plane $\BH^2\times\{0\}$. 

Now, let $\Delta^+$ be the region in the horizontal plane $\BH^2\times\{+k\}$ such that $\partial \Delta^+= \xi_0^+\cup\sigma^+_+\cup\xi_1^+\cup\sigma^+_-$ (See Figure \ref{banana}-right). Let $\Delta^-$ be the reflection of $\Delta^+$ with respect to horizontal plane $\BH^2\times\{0\}$. Let $\Omega^+$ be the region in the vertical plane containing $\sigma^+_+$ and $\sigma^-_+$ such that $\partial \Omega^+= \sigma^+_+\cup \nu_0^+\cup \sigma^-_+\cup \nu_1^+$. Similarly, define $\Omega^-$ in the opposite side. Hence, $\A=\Delta^+\cup\Delta^-\cup\Omega^+\cup\Omega^-$. Then, we have $\partial \A=\partial D_0\cup\partial D_1$.

Let $|.|$ represent the area. We claim that $|\A|<|D_0|+|D_1|$ for sufficiently large $l>0$ and $K>0$. First, note that $|D_i|>4k l$ as $2k$ is the height of the rectangle $D_i$, and any horizontal segment $\eta_t\cap D_i$ has length greater than $2l$ by construction. 

Consider $|\A|$. $\Delta^+$ belongs to the banana region in $\BH^2\times \{k\}$ bounded by $\eta_{k}$ and its reflection. Let $\beta(t)$ be the asymptotic angle between the geodesic $\tau$ and the equidistant line $\eta_t$. Note that there is a one to one correspondence between the equidistance $\rho(t)$ and the angle $\beta(t)$. Let $\beta_0=\beta(k)$. In this setting, if $t\to K$, then $\rho(t)\to \infty$ and $\beta(t)\to \frac{\pi}{2}$. 
Then, a simple computation shows that $|\Delta^+|= 4l.\tan{\beta_0}$. Furthermore, $|\Omega^\pm|<2k\|\sigma^+_+\|$ as $\Omega^+$ is a rectangle in the vertical plane with height $2k$ and all horizontal segments has length $2\rho(t)$ for $t\in[0,k]$. As $\|\sigma^+_+\|=2\rho(k)$, we have $|\Omega^\pm|<4k\rho(k)$.

Hence, we have $|\A|=2|\Delta|+2|\Omega|< 4l\tan{\beta_0} + 8k\rho(k)$

\vspace{.2cm}

Since $|D_i|>4k l$, $|\A|<|D_0|+|D_1|$ is equivalent to say that 
$$8k\rho(k)< 4l.(2k-\tan{\beta_0})$$

Now, fix $k>\frac{\pi}{2}$. Notice that by the explicit description of $\p_i$ in \cite{ST}, if the height of $\p_i$, $K\to \infty$ then $\beta_0\to 0$ and $\rho(k)\to 0$. Hence,
by choosing $K$ sufficiently large, we can make sure that $\tan{\beta_0}<2k$. Then, for sufficiently large $l>0$, we have the desired inequality. The proof of the Claim and Case 2 follows. \hfill $\Box$

\vspace{.2cm}

Now, we  finish the proof of the lemma. So far, we have shown that if $\Sigma$ is an area minimizing surface in $\BHH$ with $\PI\Sigma \subset \bigcup_{i=1}^N\gamma_{c_i}\bigcup_{j=1}^M \alpha_j$, then $\PI \Sigma \subset \bigcup_{i=1}^N\gamma_{c_i}$. In other words, we prove that the asymptotic boundary of such an area minimizing surface cannot have any vertical  segments. Now, we  show that every component of $\Sigma$ is a horizontal plane. In particular, assume that a component $S$ of $\Sigma$ contains more than 1 horizontal circle, say $\gamma_{c_1}\cup\gamma_{c_2}$. By assumption, $|c_1-c_2|>K_0>\pi$. Let $[d_1,d_2]\subset (c_1,c_2)$ with $d_2-d_1=\pi$. Then, consider the parabolic catenoid $\mathfrak{C}$ with $\PI \mathfrak{\C} = \gamma_{d_1}\cup\gamma_{d_2}\cup \alpha$ where $\alpha$ is the vertical segment corresponding to $\{0\}\times [d_1,d_2]$ in upper half space model. In particular, in the upper half space model, $\mathfrak{C}=\sigma\times\BR$ where $\sigma$ is the generating curve in $xy$-plane $\BH^2$ with $\PI \sigma = \{(d_1,0),(d_2,0)\}$. Let $\varphi_\lambda(x,y,z)=(\lambda x,\lambda y, z)$ be the isometry of $\BHH$ in the upper half space model. Then, let $\mathfrak{C}_\lambda=\varphi_\lambda(\mathfrak{C})$ is another parabolic catenoid with generating curve $\lambda\cdot \sigma$. Now, for sufficiently small $\lambda>0$, $\mathfrak{C}\cap S=\emptyset$. On the other hand, when $\lambda \to \infty$, $\mathfrak{C}$ converges to $\BH^2\times\{d_1,d_2\}$. This means if $\PI S\supset \gamma_{c_1}\cup\gamma_{c_2}$, by increasing $\lambda$, for some $\lambda_0>0$, we can find the first point of touch between $S$ and $\mathfrak{C}_{\lambda_0}$. However, this contradicts to the maximum principle. 

Finally, we  show that if $\Gamma=\bigcup_{i=1}^N\gamma_{c_i}$ ($M=0$), then $\Sigma$ is indeed a collection of horizontal planes. Assume that there is a component $S$ in $\Sigma$ with $\PI S= \gamma_{c_{j}}\cup \gamma_{c_{k}}$. Since $h(\Gamma)=K_0>\pi$, let $[e^-,e^+]\subset (c_{j},c_{k})$ with $e^+-e^-=\pi$. Let $\C$ be Daniel's parabolic catenoid with $\PI \C =\gamma_{e^+}\cup\gamma_{e^-}$. We can push $\C$ towards $\Si$ as much as we want by using isometries so that we can assume $\C\cap S=\emptyset$. Then, by pushing $\C$ towards $S$ by using the isometries, we  get a first point of touch, which contradicts to the maximum principle. This proves that $\Sigma$ must be a collection of horizontal planes, i.e. $\Sigma=\bigcup_{i=1}^N\BH^2\times\{c_i\}$. The proof follows.
\end{pf}

	
\begin{rmk} [Bridge height $K_0$] Note that the above lemma is the only reason we need large $K_0$ for the vertical bridge principle. However, the constant $K_0$ in the lemma above might be highly improved (conjecturally $K_0=\pi$) by using similar ideas. In particular, the estimates we use in Lemma 7.6 are very rough, and by using the explicit description of the generating curve for $\p_i$ in \cite{ST}, one can choose $k\in (\pi,h(\Gamma))$ more elegantly. Then, by choosing $l$ sufficiently large, one might get vertical bridge principle for all tall curves ($h(\Gamma)>\pi$), not just curves with $h(\Gamma)>K_0$. Furthermore, it might also be possible to prove a similar result for any collection of arcs $\{\alpha_i\}$ without the verticality condition on $\alpha_i$.
\end{rmk}
	



\end{document}